\documentclass[a4paper,leqno,twoside]{article}
\usepackage{latexsym,amssymb,amsmath, amsthm}
\usepackage[affil-it]{authblk}

\newtheorem{theo}{Theorem}[section]
\newtheorem{coro}[theo]{Corollary}
\newtheorem{lemm}[theo]{Lemma}
\newtheorem{prop}[theo]{Proposition}

\newtheorem{exam} [theo]{Example}
\newtheorem{rema}[theo]{Remark}

\author{I.D. Chipchakov}
\affil{Institute of Mathematics and Informatics\\ Bulgarian Academy
of Sciences\\ Acad. G. Bonchev Str., bl. 8\\ 1113, Sofia, Bulgaria;
email: chipchak@math.bas.bg}

\begin{document}
\title{On Brauer $p$-dimensions and index-exponent relations over
finitely-generated field extensions\footnote{Throughout this paper,
we write for brevity "FG-extension(s)" instead of
"finitely-generated [field] extension(s)".}}

\date{\today}
\maketitle

\begin{abstract}
Let $E$ be a field of absolute Brauer dimension abrd$(E)$, and $F/E$
a transcendental finitely-generated extension. This paper shows that
the Brauer dimension Brd$(F)$ is infinite, if abrd$(E) = \infty $.
When the absolute Brauer $p$-dimension abrd$_{p}(E)$ is infinite,
for some prime number $p$, it proves that for each pair $(n, m)$ of
integers with $n \ge m > 0$, there is a central division $F$-algebra
of Schur index $p ^{n}$ and exponent $p ^{m}$. Lower bounds on the
Brauer $p$-dimension Brd$_{p}(F)$ are obtained in some important
special cases where abrd$_{p}(E) < \infty $. These results solve
negatively a problem posed by Auel et al. (Transf. Groups {\bf 16}: 
219-264, 2011).
\end{abstract}

(secondary).

\par
\section{\bf Introduction}
\par
\medskip
Let $E$ be a field, $s(E)$ the class of finite-dimensional
associative central simple $E$-algebras, $d(E)$ the subclass of
division algebras $D \in s(E)$, and for each $A \in s(E)$, let $[A]$
be the equivalence class of $A$ in the Brauer group Br$(E)$. It is
known that Br$(E)$ is an abelian torsion group (cf. \cite{P}, Sect.
14.4), whence it decomposes into the direct sum of its
$p$-components Br$(E) _{p}$, where $p$ runs across the set $\mathbb
P$ of prime numbers. By Wedderburn's structure theorem (see, e.g.,
\cite{P}, Sect. 3.5), each $A \in s(E)$ is isomorphic to the full
matrix ring $M _{n}(D _{A})$ of order $n$ over some $D _{A} \in
d(E)$ that is uniquely determined by $A$, up-to an $E$-isomorphism.
This implies the dimension $[A\colon E]$ is a square of a positive integer deg$(A)$, the degree of $A$. The main numerical invariants
of $A$ are deg$(A)$, the Schur index ind$(A) = {\rm deg}(D _{A})$,
and the exponent exp$(A)$, i.e. the order of $[A]$ in Br$(E)$. The following statements describe basic divisibility relations between ind$(A)$ and exp$(A)$, and give an idea of their behaviour under the scalar extension map Br$(E) \to {\rm Br}(R)$, in case $R/E$ is a
field extension of finite degree $[R\colon E]$ (see, e.g., \cite{P}, Sects. 13.4, 14.4 and 15.2, and \cite{Ch2}, Lemma~3.5):
\par
\medskip\noindent
(1.1) (a) $({\rm ind}(A), {\rm exp}(A))$ is a Brauer pair, i.e.
exp$(A)$ divides ind$(A)$ and is divisible by every $p \in \mathbb
P$ dividing ind$(A)$.
\par
(b) ind$(A \otimes _{E} B)$ is divisible by ${\rm l.c.m.}\{{\rm
ind}(A), {\rm ind}(B)\}/{\rm g.c.d.}\{{\rm ind}(A), {\rm ind}(B)\}$
and divides ind$(A){\rm ind}(B)$, for each $B \in s(E)$; in
particular, if $A, B \in d(E)$ and g.c.d.$\{{\rm ind}(A), {\rm
ind}(B)\} = 1$, then the tensor product $A \otimes _{E} B$ lies in
$d(E)$.
\par
(c) ind$(A)$, ind$(A \otimes _{E} R)$, exp$(A)$ and exp$(A \otimes
_{E} R)$ divide ind$(A \otimes _{E} R)[R\colon E]$, ind$(A)$, exp$(A
\otimes _{E} R)[R\colon E]$ and exp$(A)$, respectively.
\par
\medskip
Statements (1.1) (a), (b) imply Brauer's Primary Tensor Product
Decomposition Theorem, for any $\Delta \in d(E)$ (cf. \cite{P},
Sect. 14.4), and (1.1) (a) fully describes general restrictions on
index-exponent relations, in the following sense:
\par
\medskip\noindent
(1.2) Given a Brauer pair $(m ^{\prime }, m) \in \mathbb N ^{2}$,
there is a field $F$ with $({\rm ind}(D), {\rm exp}(D)) = (m
^{\prime }, m)$, for some $D \in d(F)$ (Brauer, see \cite{P}, Sect.
19.6). One may take as $F$ any rational (i.e. purely transcendental)
extension in infinitely many variables over any fixed field $F _{0}$
(see also Corollary \ref{coro4.4} and Remark \ref{rema4.5}).
\par
\medskip
As in \cite{ABGV}, Sect. 4, we say that a field $E$ is of finite
Brauer $p$-dimension Brd$_{p}(E) = n$, for a fixed $p \in \mathbb
P$, if $n$ is the least integer $\ge 0$, for which ind$(D) \le {\rm
exp}(D) ^{n}$ whenever $D \in d(E)$ and $[D] \in {\rm Br}(E) _{p}$.
If no such $n$ exists, we set Brd$_{p}(E) = \infty $. The absolute
Brauer $p$-dimension of $E$ is defined as the supremum abrd$_{p}(E)
= {\rm sup}\{{\rm Brd}_{p}(R)\colon \ R \in {\rm Fe}(E)\}$, where
Fe$(E)$ is the set of finite extensions of $E$ in a separable
closure $E _{\rm sep}$. Clearly, Brd$_{p}(E) \le {\rm abrd}_{p}(E)$,
$p \in \mathbb P$. We say that $E$ is a virtually perfect field, if char$(E) = 0$ or char$(E) = q > 0$ and $E$ is a finite
extension of its subfield $E ^{q} = \{e ^{q}\colon \ e \in E\}$.
\par
\medskip
It is known that ${\rm Brd}_{p}(E) = {\rm abrd}_{p}(E) = 1$, for all
$p \in \mathbb P$, if $E$ is a global or local field (cf. \cite{Re},
(31.4) and (32.19)), or the function field of an algebraic surface
defined over an algebraically closed field $E _{0}$ \cite{Jong},
\cite{Lieb} (see also Remark \ref{rema5.8}). As shown in \cite{Mat},
abrd$_{p}(E) < p ^{n-1}$, $p \in \mathbb P$, provided that $E$ is the function field of an $n$-dimensional algebraic variety defined over
an algebraically closed field $E _{0}$. Similarly, abrd$_{p}(E) < p ^{n}$, $p \in \mathbb P$, if $E _{0}$ is a finite field, the maximal unramified extension of a local field, or a perfect pseudo
algebraically closed (PAC) field (for the $C _{1}$-type of $E _{0}$,
used in \cite{Mat} for proving these inequalities, see \cite{L1} and \cite{Kol}, \cite{FJ}, Theorem~21.3.6, respectively). The suprema
Brd$(E) = {\rm sup}\{{\rm Brd}_{p}(E)\colon \ p \in \mathbb P\}$ and abrd$(E) = {\rm sup}\{{\rm Brd}(R)\colon \ R \in
{\rm Fe}(E)\}$ are called a Brauer dimension and an absolute Brauer dimension of $E$, respectively. In view of (1.1), the definition of Brd$(E)$ is the same as the one given in \cite{ABGV}, Sect. 4. It has recently been proved \cite{HHKr}, \cite{PaSur} (see also \cite{Ch5},
Propositions~6.1 and 7.1), that abrd$(K _{m}) < \infty $, provided
$m \in \mathbb N$ and $(K _{m}, v _{m})$ is an $m$-dimensional local
field, in the sense of \cite{FV}, with a finite $m$-th residue field
$\widehat K _{m}$.
\par
\medskip
The present research is devoted to the study of index-exponent
relations over transcendental FG-extensions $F$ of a field $E$ and
their dependence on abrd$_{p}(E)$, $p \in \mathbb P$. It is
motivated mainly by two questions concerning the dependence of
Brd$(F)$ upon Brd$(E)$, stated as open problems in Sect. 4 of the
survey \cite{ABGV}.

\medskip
\section{\bf The main results}
\par
\medskip
While the study of index-exponent relations makes interest in its
own right, it should be noted that fields $E$ with abrd$_{p}(E) <
\infty $, for all $p \in \mathbb P$, are singled out by Galois
cohomology (see \cite{Ka} and \cite{Vo}, as well as \cite{Mat},
Sects. 5-8, and further references in \cite{Ch4}, Remark~4.2). It is
also worth mentioning the following fact about the almost perfect
fields of this type (see \cite{Ch1}, \cite{Ch2}, and Lemma \ref{lemm4.1}):
\par
\medskip
(2.1) Every locally finite dimensional associative central division
$E$-algebra $R$ possesses an $E$-subalgebra $\widetilde R$ with the
following properties:
\par
(a) $\widetilde R$ decomposes into a tensor product $\otimes _{p \in
\mathbb P} R _{p}$, where $\otimes = \otimes _{E}$, $R _{p} \in
d(E)$ and $[R _{p}] \in {\rm Br}(E) _{p}$, for each $p \in \mathbb
P$;
\par
(b) Finite-dimensional $E$-subalgebras of $R$ are embeddable in
$\widetilde R$;
\par
(c) $\widetilde R$ is isomorphic to $R$, if the dimension $[R\colon
E]$ is countably infinite.
\par
\medskip
It would be of definite interest to know whether function fields of
algebraic varieties over a global, local or algebraically closed
field are of finite absolute Brauer dimensions. This draws our
attention to the following open question:
\par
\medskip\noindent
(2.2) Is the class of fields $E$ of finite absolute Brauer
$p$-dimensions, for a fixed $p \in \mathbb P$, $p \neq {\rm
char}(E)$, closed under the formation of FG-extensions?
\par
\medskip
The main result of this paper shows, for a transcendental
FG-extension $F/E$, the strong influence of $p$-dimensions
abrd$_{p}(E)$ on Brd$_{p}(F)$, and on index-exponent relations over
$F$, as follows:
\par
\medskip
\begin{theo}
\label{theo2.1} Let $E$ be a field, $p \in \mathbb P$ and $F/E$ an
{\rm FG}-extension of transcendency degree {\rm trd}$(F/E) = \kappa
\ge 1$. Then:
\par
{\rm (a)} {\rm Brd}$_{p}(F) \ge {\rm abrd}_{p}(E) + \kappa - 1$, if
{\rm abrd}$_{p}(E) < \infty $ and $F/E$ is rational;
\par
{\rm (b)} If {\rm abrd}$_{p}(E) = \infty $, then {\rm Brd}$_{p}(F) =
\infty $ and for each $n, m \in \mathbb N$ with $n \ge m > 0$, there
exists $D _{n,m} \in d(F)$ with {\rm ind}$(D _{n,m}) = p ^{n}$ and
{\rm exp}$(D _{n,m}) = p ^{m}$;
\par
{\rm (c)} Brd$_{p}(F) = \infty $, provided $p = {\rm char}(E)$ and
$[E\colon E ^{p}] = \infty $; if {\rm char}$(E) = p$ and $[E\colon E
^{p}] = p ^{\nu } < \infty $, then $\nu + \kappa - 1 \le {\rm
Brd}_{p}(F) \le {\rm abrd}_{p}(F) \le \nu + \kappa $.
\end{theo}
\par
\medskip
It is known (cf. \cite{L}, Ch. X) that each FG-extension $F$ of a
field $E$ possesses a subfield $F _{0}$ that is rational over $E$
with trd$(F _{0}/E) = {\rm trd}(F/E)$. This ensures that $[F\colon F
_{0}] < \infty $, so (1.1) and Theorem \ref{theo2.1} imply the
following:
\par
\medskip\noindent
(2.3) If (2.2) has an affirmative answer, for some $p \in \mathbb
P$, $p \neq {\rm char}(E)$, and each FG-extension $F/E$ with
trd$(F/E) = \kappa \ge 1$, then there exists $c _{\kappa }(p) \in
\mathbb N$, depending on $E$, such that Brd$_{p}(\Phi ) \le c
_{\kappa }(p)$, for every FG-extension $\Phi /E$ with trd$(\Phi
/E) < \kappa $. For example, this applies to $c _{k}(p) = {\rm
Brd}_{p}(E _{\kappa })$, where $E _{\kappa }/E$ is a rational
FG-extension with trd$(E _{\kappa }/E) = \kappa $.
\par
\vskip0.34truecm
The application of Theorem \ref{theo2.1} is facilitated by the
following result of \cite{Ch4} (see Example \ref{exam6.2} below,
for an alternative proof in characteristic zero):
\par
\medskip
\begin{prop}
\label{prop2.2} For each $q \in \mathbb P \cup \{0\}$ and $k \in
\mathbb N$, there exists a field $E _{q,k}$ with char$(E _{q,k}) =
q$, Brd$(E _{q,k}) = k$ and abrd$_{p}(E _{q,k}) = \infty $, for all
$p \in \mathbb P \setminus P _{q}$, where $P _{0} = \{2\}$ and $P
_{q} = \{p \in \mathbb P\colon \ p \mid q(q - 1)\}$, $q \in \mathbb
P$. Moreover, if $q > 0$, then $E _{q,k}$ can be chosen so that $[E
_{q,k}\colon E _{q,k} ^{q}] = \infty $.
\end{prop}
\par
\medskip
Theorem \ref{theo2.1}, Proposition \ref{prop2.2} and statement (1.1)
(b) imply the following:
\par
\medskip\noindent
(2.4) There exist fields $E _{k}$, $k \in \mathbb N$, such that
char$(E _{k}) = 2$, Brd$(E _{k}) = k$ and all Brauer pairs $(m
^{\prime }, n ^{\prime }) \in \mathbb N ^{2}$ are index-exponent
pairs over any transcendental FG-extension of $E _{k}$.
\par
\medskip\noindent
It is not known whether (2.4) holds in any characteristic $q \neq
2$. This is closely related to the following open problem:
\par
\medskip\noindent
(2.5) Find whether there exists a field $E$ containing a primitive
$p$-th root of unity, for a given $p \in \mathbb P$, such that
Brd$_{p}(E) < {\rm abrd}_{p}(E) = \infty $.
\par
\medskip
Statement (1.1) (b), Theorem \ref{theo2.1} and Proposition
\ref{prop2.2} imply the validity of (2.4) in zero characteristic,
for Brauer pairs of odd positive integers. When $q > 2$, they show
that if $[E _{q,k}\colon E _{q,k} ^{q}] = \infty $, then Brauer
pairs $(m ^{\prime }, m) \in \mathbb N ^{2}$ relatively prime to $q
- 1$ are index-exponent pairs over every transcendental
FG-extension of $E _{q,k}$. This solves in the negative \cite{ABGV},
Problem~4.4, proving (in the strongest presently known form) that
the class of fields of finite Brauer dimensions is not closed under
the formation of FG-extensions.
\par
\medskip
Theorem \ref{theo2.1} (a) makes it easy to prove that the solution
to \cite{ABGV}, Problem~4.5, on the existence of a "good" definition
of a dimension dim$(E) < \infty $, for some fields $E$, is negative
whenever abrd$(E) = \infty $ (see Corollary \ref{coro5.4}). It
implies that if Problem~4.5 of \cite{ABGV} is solved affirmatively,
for all FG-extensions $F/E$, then each $F$ satisfies the following 
stronger inequalities than those conjectured by (2.3) (see also 
Remark \ref{rema5.5} and \cite{ABGV}, Sect. 4):
\par
\medskip\noindent
(2.6) Brd$(F) < {\rm dim}(F)$, abrd$(F) \le {\rm dim}(F)$ and
abrd$(F) \le {\rm Brd}(E _{t+1}) \le {\rm abrd}(E) + t + c(E)$, for
some integer $c(E) \le {\rm dim}(E) - {\rm abrd}(E)$, where $t =
{\rm trd}(F/E)$, $E _{t+1}/E$ is a rational extension and trd$(E
_{t+1}/E) = t + 1$.
\par
\medskip
The proof of Theorem \ref{theo2.1} is based on Merkur'ev's theorem
about central division algebras of prime exponent \cite{M1}, Sect.
4, Theorem~2, and on a characterization of fields of finite absolute
Brauer $p$-dimensions generalizing Albert's theorem \cite{A1}, Ch.
XI, Theorem~3. It strongly relies on results of valuation theory,
like theorems of Grunwald-Hasse-Wang type, Morandi's theorem on
tensor products of valued division algebras \cite{Mo}, Theorem~1,
lifting theorems over Henselian (valued) fields, and Ostrowski's
theorem. As shown in \cite{Ch4}, Sect. 6, the flexibility of this
approach enables one to obtain the following results:
\par
\medskip\noindent
(2.7) (a) There exists a field $E _{1}$ with abrd$(E _{1}) = \infty
$, abrd$_{p}(E _{1}) < \infty $, $p \in \mathbb P$, and Brd$(L _{1})
< \infty $, for every finite extension $L _{1}/E _{1}$;
\par
(b) For any integer $n \ge 2$, there is a Galois extension $L _{n}/E
_{n}$, such that $[L _{n}\colon E _{n}] = n$, Brd$_{p}(L _{n}) =
\infty $, for all $p \in \mathbb P$, $p \equiv 1 ({\rm mod} \ n)$,
and Brd$(M _{n}) < \infty $, provided that $M _{n}$ is an extension
of $E$ in $L _{n,{\rm sep}}$ not including $L _{n}$.
\par
\medskip
Our basic notation and terminology are standard. For any field $K$ 
with a Krull valuation $v$, unless 
stated otherwise, we denote by $O _{v}(K)$, $\widehat K$ and $v(K)$
the valuation ring, the residue field and the value group of $(K,
v)$, respectively; $v(K)$ is supposed to be an additively written
totally ordered abelian group. As usual, $\mathbb Z$ stands for the
additive group of integers, $\mathbb Z _{p}$, $p \in \mathbb P$, are
the additive groups of $p$-adic integers, and $[r]$ is the integral
part of any real number $r \ge 0$. We write $I(\Lambda ^{\prime
}/\Lambda )$ for the set of intermediate fields of a field extension
$\Lambda ^{\prime }/\Lambda $, and Br$(\Lambda ^{\prime }/\Lambda )$
for the relative Brauer group of $\Lambda ^{\prime }/\Lambda $. By a
$\Lambda $-valuation of $\Lambda ^{\prime }$, we mean a Krull
valuation $v$ with $v(\lambda ) = 0$, for all $\lambda \in
\Lambda ^{\ast }$. Given a field $E$ and $p \in \mathbb P$, $E(p)$
denotes the maximal $p$-extension of $E$ in $E _{\rm sep}$, and $r
_{p}(E)$ the rank of the Galois group $\mathcal{G}(E(p)/E)$ as a
pro-$p$-group ($r _{p}(E) = 0$, if $E(p) = E$). Brauer groups are
considered to be additively written, Galois groups are viewed as
profinite with respect to the Krull topology, and by a homomorphism
of profinite groups, we mean a continuous one. We refer the reader
to \cite{Efr2}, \cite{JW}, \cite{L}, \cite{P} and \cite{Se}, for any
missing definitions concerning valuation theory, field extensions,
simple algebras, Brauer groups and Galois cohomology.
\par
\medskip
The rest of the paper proceeds as follows: Sect. 3 includes
preliminaries used in the sequel. Theorem \ref{theo2.1} is proved in
Sects. 4 and 5. In Sect. 6 we show that the answer to (2.2) will 
be affirmative, if this is the case in zero characteristic. Lower 
bounds on Brd$_{p}(F)$ are also obtained in these sections, for 
FG-extensions $F$ of some frequently used fields $E$ with 
abrd$_{p}(E) < \infty $.

\medskip
\section{\bf Preliminaries on valuation theory}
\par
\medskip
The results of this section are known and will often be used without
an explicit reference. We begin with a lemma essentially due to
Saltman \cite{Sal2}.
\par
\medskip
\begin{lemm}
\label{lemm3.1} Let $(K, v)$ be a height $1$ valued field, $K _{v}$
a Henselization of $K$ in $K _{\rm sep}$ relative to $v$, and
$\Delta _{v} \in d(K _{v})$ an algebra of exponent $p \in \mathbb
P$. Then there exists $\Delta \in d(K)$ with {\rm exp}$(\Delta ) =
p$ and $[\Delta \otimes _{K} K _{v}] = [\Delta _{v}]$.
\end{lemm}

\medskip
\begin{proof}
By \cite{M1}, Sect. 4, Theorem~2, $\Delta _{v}$ is Brauer equivalent
to a tensor product of degree $p$ algebras from $d(K _{v})$, so one
may consider only the case of deg$(\Delta _{v}) = p$. Then, by
Saltman's theorem (cf. \cite{Sal2}), there exists $\Delta \in d(K)$,
such that deg$(\Delta ) = p$ and $\Delta \otimes _{K} K _{v}$ is $K
_{v}$-isomorphic to $\Delta _{v}$, which proves Lemma \ref{lemm3.1}.
\end{proof}

\medskip
In what follows, we shall use the fact that the Henselization $K
_{v}$ of a field $K$ with a valuation $v$ of height $1$ is separably
closed in the completion of $K$ relative to the topology induced by
$v$ (cf. \cite{Efr2}, Theorem~15.3.5 and Sect. 18.3). For example,
our next lemma is a consequence of Galois theory, this fact and
Lorenz-Roquette's valuation-theoretic generalization of
Grunwald-Wang's theorem (cf. \cite{L}, Ch. VIII, Theorem~4, and
\cite{LR}, page 176 and Theorems~1 and 2).

\medskip
\begin{lemm}
\label{lemm3.2}
Let $F$ be a field, $S = \{v _{1}, \dots , v _{s}\}$ a finite set of
non-equivalent height $1$ valuations of $F$, and for each index $j$,
let $F _{v _{j}}$ be a Henselization of $K$ in $K _{\rm sep}$
relative to $v _{j}$, and $L _{j}/F _{v _{j}}$ a cyclic field
extension of degree $p ^{\mu _{j}}$, for some $p \in P$
and $\mu _{j} \in \mathbb N$. Put $\mu = {\rm max}\{\mu _{1}, \dots
, \mu _{s}\}$, and in the case of $p = 2$ and {\rm char}$(F) = 0$,
suppose that the extension $F(\delta _{\mu })/F$ is cyclic, where
$\delta _{\mu } \in F _{\rm sep}$ is a primitive $2 ^{\mu }$-th root
of unity. Then there is a cyclic field extension $L/F$ of degree $p
^{\mu }$, whose Henselization $L _{v _{j}'}$ is $F _{v
_{j}}$-isomorphic to $L _{j}$, where $v _{j} ^{\prime }$ is a
valuation of $L$ extending $v _{j}$, for $j = 1, \dots , s$.
\end{lemm}

\medskip
Assume that $K = K _{v}$, or equivalently, that $(K, v)$ is a
Henselian field, i.e. $v$ is a Krull valuation on $K$, which extends
uniquely, up-to an equivalence, to a valuation $v _{L}$ on each
algebraic extension $L/K$. Put $v(L) = v _{L}(L)$ and denote by
$\widehat L$ the residue field of $(L, v _{L})$. It is known that
$\widehat L/\widehat K$ is an algebraic extension and $v(K)$ is a
subgroup of $v(L)$. When $[L\colon K]$ is finite, Ostrowski's
theorem states the following (cf. \cite{Efr2}, Theorem~17.2.1):
\par
\medskip\noindent
(3.1) $[\widehat L\colon \widehat K]e(L/K)$ divides $[L\colon K]$ and
$[L\colon K][\widehat L\colon \widehat K] ^{-1}e(L/K) ^{-1}$ is not
divisible by any $p \in \mathbb P$ different from char$(\widehat K)$,
$e(L/K)$ being the index of $v(K)$ in $v(L)$;  in particular, if
char$(\widehat K) \dagger [L\colon K]$, then $[L\colon K] = [\widehat
L\colon \widehat K]e(L/K)$.
\par
\medskip
Statement (3.1) and the Henselity of $v$ imply the following:
\par
\medskip\noindent
(3.2) The quotient groups $v(K)/pv(K)$ and $v(L)/pv(L)$ are
isomorphic, if $p \in \mathbb P$ and $L/K$ is a finite extension.
When char$(\widehat K) \dagger [L\colon K]$, the natural embedding
of $K$ into $L$ induces canonically an isomorphism $v(K)/pv(K) \cong
v(L)/pv(L)$.
\par
\medskip
A finite extension $R/K$ is said to be defectless, if $[R\colon K] =
[\widehat R\colon \widehat K]e(R/K)$. It is called inertial, if
$[R\colon K] = [\widehat R\colon \widehat K]$ and $\widehat R$ is
separable over $\widehat K$. We say that $R/K$ is totally ramified,
if $[R\colon K] = e(R/K)$; $R/K$ is called tamely ramified, if
$\widehat R/\widehat K$ is separable and char$(\widehat K) \dagger
e(R/K)$. The Henselity of $v$ ensures that the compositum $K _{\rm
ur}$ of inertial extensions of $K$ in $K _{\rm sep}$ has the
following properties:
\par
\medskip\noindent
(3.3) (a) $v(K _{\rm ur}) = v(K)$ and finite extensions of $K$ in $K
_{\rm ur}$ are inertial;
\par
(b) $K _{\rm ur}/K$ is a Galois extension, $\widehat K _{\rm ur}
\cong \widehat K _{\rm sep}$ over $\widehat K$, $\mathcal{G}(K _{\rm
ur}/K)$ $\cong \mathcal{G}_{\widehat K}$, and the natural mapping of
$I(K _{\rm ur}/K)$ into $I(\widehat K _{\rm sep}/\widehat K)$ is
bijective.
\par
\medskip\noindent
Recall that the compositum $K _{\rm tr}$ of tamely ramified
extensions of $K$ in $K _{\rm sep}$ is a Galois extension of $K$
with $v(K _{\rm tr}) = pv(K _{\rm tr})$, for every $p \in \mathbb P$
not equal to char$(\widehat K)$. It is therefore clear from (3.1)
that if $K _{\rm tr} \neq K _{\rm sep}$, then char$(\widehat K) = q
\neq 0$ and $\mathcal{G}_{K _{\rm tr}}$ is a pro-$q$-group. When
this holds, it follows from (3.3) and Galois cohomology (cf.
\cite{Se}, Ch. II, 2.2) that cd$_{q}(\mathcal{G}(K _{\rm tr}/K)) \le
1$. Hence, by \cite{Se}, Ch. I, Proposition~16, there is a closed
subgroup $\mathcal{H} \le \mathcal{G}_{K}$, such that
$\mathcal{G}_{K _{\rm tr}}\mathcal{H} = \mathcal{G}_{K}$,
$\mathcal{G}_{K _{\rm tr}} \cap \mathcal{H} = \{1\}$ and
$\mathcal{H} \cong \mathcal{G}(K _{\rm tr}/K)$. In view of Galois
theory and the Mel'nikov-Tavgen' theorem \cite{MT}, these results
imply in the case of char$(\widehat K) = q > 0$ the existence of a
field $K ^{\prime } \in I(K _{\rm sep}/K)$ satisfying the following
conditions:
\par
\medskip\noindent
(3.4) $K ^{\prime } \cap K _{\rm tr} = K$, $K ^{\prime }K _{\rm tr}
= K _{\rm sep}$ and $K _{\rm sep} \cong K _{\rm tr} \otimes _{K} K
^{\prime }$ over $K$; the field $\widehat K ^{\prime }$ is a perfect
closure of $\widehat K$, finite extensions of $K$ in $K ^{\prime }$
are of $q$-primary degrees, $K _{\rm sep} = K ^{\prime }_{\rm tr}$,
$v(K ^{\prime }) = qv(K ^{\prime })$, and the natural embedding of
$K$ into $K ^{\prime }$ induces isomorphisms $v(K)/pv(K) \cong v(K
^{\prime })/pv(K ^{\prime })$, $p \in \mathbb P \setminus \{q\}$.

\medskip
Assume as above that $(K, v)$ is Henselian. Then each $\Delta \in
d(K)$ has a unique, up-to an equivalence, valuation $v _{\Delta }$
extending $v$ so that the value group $v(\Delta )$ of $(\Delta , v
_{\Delta })$ is totally ordered and abelian (cf. \cite{Sch}, Ch. 2,
Sect. 7). It is known that $v(K)$ is a subgroup of $v(\Delta )$ of
index $e(\Delta /K) \le [\Delta \colon K]$, and the residue division
ring $\widehat {\Delta }$ of $(\Delta , v _{\Delta })$ is a
$\widehat K$-algebra. Moreover, by the Ostrowski-Draxl theorem
\cite{Dr2}, $[\Delta \colon K]$ is divisible by $e(\Delta
/K)[\widehat {\Delta }\colon \widehat K]$, and in case
char$(\widehat K) \dagger [\Delta \colon K]$, $[\Delta \colon K] =
e(\Delta /K)[\widehat \Delta \colon \widehat K]$. An algebra $D \in
d(K)$ is called inertial, if $[D\colon K] = [\widehat D\colon
\widehat K]$ and $\widehat D \in d(\widehat K)$. Inertial
$K$-algebras and algebras from $d(\widehat K)$ are related as follows  (see \cite{JW}, Theorem~2.8):
\par
\medskip\noindent
(3.5) (a) Each $\widetilde D \in d(\widehat K)$ has an inertial lift
over $K$, i.e. $\widetilde D = \widehat D$, for some $D \in d(K)$
inertial over $K$, that is uniquely determined by $\widetilde D$,
up-to a $K$-isomorphism.
\par
(b) The set IBr$(K) = \{[I] \in {\rm Br}(K)\colon \ I \in d(K)$ is
inertial$\}$ is a subgroup of Br$(K)$; the canonical map IBr$(K)
\to {\rm Br}(\widehat K)$ is an index-preserving isomorphism.

\medskip
\section{\bf Proof of Theorem \ref{theo2.1} (a) and (c)}

\medskip
The role of Lemma \ref{lemm3.1} in the study of Brauer
$p$-dimensions of FG-extensions of a field $E$ is determined by the
following result of \cite{Ch4}, which characterizes the condition
abrd$_{p}(E) \le \mu $, for a given $\mu \in \mathbb N$. When $E$ is
virtually perfect, this result is in fact equivalent to
\cite{PaSur}, Lemma~1.1, and in case $\mu = 1$, it restates
Theorem~3 of \cite{A1}, Ch. XI.

\medskip
\begin{lemm}
\label{lemm4.1} Let $E$ be a field, $p \in \mathbb P$ and $\mu \in
\mathbb N$. Then {\rm abrd}$_{p}(E) \le \mu $ if and only if, for
each $E ^{\prime } \in {\rm Fe}(E)$, {\rm ind}$(\Delta ) \le p ^{\mu
}$ whenever $\Delta \in d(E ^{\prime })$ and {\rm exp}$(\Delta ) =
p$. Moreover, if $E$ is virtually perfect, then {\rm abrd}$_{p}(E)
\ge {\rm Brd}_{p}(E ^{\prime })$, for all finite extensions $E
^{\prime }/E$.
\end{lemm}

\medskip
Let now $F/E$ be a transcendental FG-extension and $F _{0} \in
I(F/E)$ a rational extension of $E$ with trd$(F _{0}/E) = {\rm
trd}(F/E) = t$. Clearly, an ordering on a fixed transcendency basis
of $F _{0}/E$ gives rise to a height $t$ $E$-valuation $v _{0}$ of $F
_{0}$ with $v _{0}(F _{0}) = \mathbb Z ^{t}$ and $\widehat F
_{0} = E$. Considering any prolongation of $v _{0}$ on $F$, and
taking into account that $[F\colon F _{0}] < \infty $, one obtains
the following:
\par
\medskip\noindent
(4.1) $F$ has an $E$-valuation $v$ of height $t$, such that $v(F)
\cong \mathbb Z ^{t}$ and $\widehat F$ is a finite extension of $E$;
in particular, $v(F)/pv(F)$ is a group of order $p ^{t}$, for every
$p \in \mathbb P$.

\medskip\noindent
When char$(E) = p$, (4.1) implies $[\widehat F\colon \widehat F
^{p}] = [E\colon E ^{p}]$ (cf. \cite{L}, Ch. VII, Sect. 7), so 
the former assertion of Theorem \ref{theo2.1} (c) follows from 
the next lemma.

\medskip
\begin{lemm}
\label{lemm4.2} Let $(K, v)$ be a valued field with {\rm char}$(K) =
q > 0$ and $v(K) \neq qv(K)$, and let $\tau (q)$ be the dimension of
$v(K)/qv(K)$ as a vector space over the field $\mathbb F _{q}$ with
$q$ elements. Then:
\par
{\rm (a)} For each $\pi \in K ^{\ast }$ with $v(\pi ) \notin qv(K)$,
there are degree $q$ extensions $L _{m}$ of $K$ in $K(q)$, $m \in
\mathbb N$, such that the compositum $M _{m} = L _{1} \dots L _{m}$
has a unique valuation $v _{m}$ extending $v$, up-to an equivalence,
$(M _{m}, v _{m})/(K, v)$ is totally ramified, $[M _{m}\colon K] = q
^{m}$ and $v(\pi ) \in q ^{m}v _{m}(M _{m})$, for each $m$;
\par
{\rm (b)} Given an integer $n \ge 2$, there exists $T _{n} \in
d(K)$ with {\rm exp}$(T _{n}) = q$ and {\rm ind}$(T _{n}) = q
^{n-1}$ except, possibly, if $\tau (q) < \infty $ and $[\widehat
K\colon \widehat K ^{q}] < q ^{n-\tau (q)}$.
\end{lemm}

\medskip
\begin{proof}
It suffices to consider the special case of $v(\pi ) < 0$. Fix a
Henselization $(K _{v}, \bar v)$ of $(K, v)$, put $\rho (K _{v}) =
\{u ^{q} - u\colon \ u \in K _{v}\}$, and for each $m \in \mathbb
N$, denote by $L _{m}$ the root field in $K _{\rm sep}$ over $K$ of
the polynomial $f _{m}(X) = X ^{q} - X - \pi _{m}$, where $\pi _{m}
= \pi ^{1+qm}$. Also, let $\mathbb F$ be the prime subfield of $K$,
$\Phi = \mathbb F(\pi )$, $\omega $ the valuation of $\Phi $ induced
by $v$, and $(\Phi _{\omega }, \bar \omega )$ a Henselization of
$(\Phi , \omega )$, such that $\Phi _{\omega } \subseteq K _{v}$ and
$\bar v$ extends $\bar \omega $ (the existence of $(\Phi _{\omega },
\bar \omega )$ follows from \cite{Efr2}, Theorem~15.3.5).
Identifying $K _{v}$ with its $K$-isomorphic copy in $K _{\rm sep}$,
put $L _{m} ^{\prime } = L _{m}K _{v}$ and $M _{m} ^{\prime } = M
_{m}K _{v}$, for every index $m$. It is easily verified that $\rho
(K _{v})$ is an $\mathbb F$-subspace of $K _{v}$ and $\bar v(u ^{q}
- u) \in q\bar v(K _{v})$, for every $u \in K _{v}$ with $\bar v(u)
< 0$. As $\bar v(K _{v}) = v(K)$, this observation and the choice of
$\pi $ indicate that the cosets $\pi _{m} + \rho (K _{v})$, $m \in
\mathbb N$, are linearly independent over $\mathbb F$. In view of
the Artin-Schreier theorem and Galois theory (cf. \cite{L}, Ch.
VIII, Sect. 6), this implies $f _{m}(X)$ is irreducible over $K
_{v}$, $L _{m} ^{\prime }/K _{v}$ and $L _{m}/K$ are cyclic
extensions of degree $q$, $M _{m} ^{\prime }/K _{v}$ and $M _{m}/K$
are abelian, and $[M _{m} ^{\prime }\colon K _{v}] = [M _{m}\colon
K] = q ^{m}$, for each $m \in \mathbb N$. Moreover, our argument
proves that degree $q$ extensions of $K _{v}$ in the compositum of
the fields $L _{m} ^{\prime }$, $m \in \mathbb N$, are cyclic and
totally ramified over $K _{v}$. At the same time, it follows from
the Henselity of $\bar v$ and the equality $\widehat K _{v} =
\widehat K$ that $M _{m} ^{\prime }$ contains as a subfield an
inertial lift over $K _{v}$ of the separable closure of $\widehat K$
in $\widehat M _{m} ^{\prime }$. When $v$ is discrete and $\widehat
K$ is perfect, the obtained results imply the assertions of Lemma
\ref{lemm4.2} (a), since finite extensions of $K _{v}$ in $K _{\rm
sep}$ are defectless (relative to $\bar v$, see \cite{L}, Ch. XII,
Sect. 6, Corollary~2).
\par
To prove Lemma \ref{lemm4.2} (a) in general it remains to be seen
that, for any fixed $m \in \mathbb N$, $M _{m}$ has a unique, up-to
an equivalence, valuation $v _{m}$ extending $v$, $(M _{m}, v
_{m})/(K, v)$ is totally ramified and $v(\pi ) \in q ^{m}v(M _{m})$.
The extendability of $v$ to a valuation $v _{m}$ of $M _{m}$ is
well-known (cf. \cite{L}, Ch. XII, Sect. 4), so our assertions can
be deduced from the concluding one, the equality $[M _{m}\colon K] =
[M _{m}K _{v}\colon K _{v}] = q ^{m}$ and statement (3.1). Our proof
also relies on the fact that $(\Phi , \omega )$ is a discrete valued
field and $\widehat \Phi /\mathbb F$ is a finite extension (see
\cite{CF}, Ch. II, Lemma~3.1, or \cite{Efr2}, Example~4.1.3); in
particular, $\widehat \Phi $ is perfect. Let now $\Psi _{m} \in I(K
_{\rm sep}/\Phi )$ be the root field of $f _{m}(X)$ over $\Phi $.
Then $L _{m} = \Psi _{m}K$, $[\Psi _{m}\colon \Phi ] = q$, $M _{m} =
\Theta _{m}K$ and $[\Theta _{m}\colon \Phi ] = q ^{m}$, where
$\Theta _{m} = \Psi _{1} \dots \Psi _{m}$. Therefore, $\Theta
_{m}\Phi _{\omega }/\Phi _{\omega }$ is totally ramified relative to
$\bar \omega $. Equivalently, the integral closure of $O _{\omega
}(\Phi )$ in $\Theta _{m}$ contains a primitive element $t _{m}
^{\prime }$ of $\Theta _{m}/\Phi $, whose minimal polynomial $\theta
_{m}(X)$ over $O _{\omega }(\Phi )$ is Eisensteinian (cf. \cite{CF},
Ch. I, Theorem~6.1, and \cite{L}, Ch. XII, Sects. 2, 3 and 6).
Hence, $\omega $ has a unique prolongation $\omega _{m}$ on $\Theta
_{m}$, up-to an equivalence, $\omega (t _{m}) \notin q\omega (\Phi
)$ and $q ^{m}\omega _{m}(t _{m} ^{\prime }) = \omega (t _{m})$,
where $t _{m}$ is the free term of $\theta _{m}(X)$. As $\pi \in
\Phi $, $v(\pi ) \notin qv(K)$ and $\Theta _{m}/\Phi $ is a Galois
extension, this implies $t _{m} ^{\prime }$ is a primitive element
of $M _{m}/K$ and $M _{m} ^{\prime }/K _{v}$, $q ^{m}v _{m}(t _{m}
^{\prime }) = v(t _{m}) = \omega (t _{m})$ and $v(\pi ) \in q ^{m}v
_{m}(M _{m})$, which completes the proof of Lemma \ref{lemm4.2} (a).
\par
We prove Lemma \ref{lemm4.2} (b). Put $\pi _{1} = \pi $ and suppose
that there exist elements $\pi _{j} \in K ^{\ast }$, $j = 2, \dots ,
n$, and an integer $\mu \le n$, such that the cosets $v(\pi _{i}) +
qv(K)$, $i = 1, \dots , \mu $, are linearly independent over
$\mathbb F _{q}$, and in case $\mu < n$, $v(\pi _{u}) = 0$ and the
residue classes $\hat \pi _{u}$, $u = \mu + 1, \dots , n$, generate
an extension of $\widehat K ^{q}$ of degree $q ^{n-\mu }$. Fix a
generator $\lambda _{m}$ of $\mathcal{G}(L _{m}/K)$, for each $m \in
\mathbb N$, denote by $T _{n}$ the $K$-algebra $\otimes _{j=2} ^{n}
(L _{j-1}/K, \lambda _{j-1}, \pi _{j})$, where $\otimes = \otimes
_{K}$, and put $T _{n} ^{\prime } = T _{n} \otimes _{K} K _{v}$. We
show that $T _{n} \in d(K)$ (whence exp$(T _{n}) = q$ and ind$(T
_{n}) = q ^{n-1}$). Clearly, there is a $K _{v}$-isomorphism $T _{n}
^{\prime } \cong \otimes _{j=2} ^{n} (L _{j-1} ^{\prime }/K _{v},
\lambda _{j-1} ^{\prime }, \pi _{j})$, where $\otimes = \otimes _{K
_{v}}$ and $\lambda _{j-1} ^{\prime }$ is the unique $K
_{v}$-automorphism of $L _{j-1} ^{\prime }$ extending $\lambda
_{j-1}$, for each $j$. Therefore, it suffices for the proof of Lemma
\ref{lemm4.2} (b) to show that $T _{n} ^{\prime } \in d(K _{v})$.
Since $K _{v}$ and $L _{m} ^{\prime }$, $m \in \mathbb N$, are
related as $K$ and $L _{m}$, $m \in \mathbb N$, this amounts to
proving that $T _{n} \in d(K)$, for $(K, v)$ Henselian. Suppose
first that $n = 2$. As $L _{1}/K$ is totally ramified, it follows
from the Henselity of $v$ that $v(l) \in qv(L _{1})$, for every
element $l$ of the norm group $N(L _{1}/K)$. One also concludes that
if $l \in N(L _{1}/K)$ and $v _{L}(l) = 0$, then $\hat l \in
\widehat K ^{q}$. These observations prove that $\pi _{2} \notin N(L
_{1}/K)$, so it follows from \cite{P}, Sect. 15.1, Proposition~b,
that $T _{2} \in d(K)$. Henceforth, we assume that $n \ge 3$ and
view all value groups considered in the rest of the proof as
(ordered) subgroups of a fixed divisible hull of $v(K)$. Note that
the centralizer $C _{n}$ of $L _{n}$ in $T _{n}$ is $L
_{n}$-isomorphic to $T _{n-1} \otimes _{K} L _{n}$ and $\otimes
_{j=2} ^{n-1} (L _{j-1}L _{n}, \lambda _{j-1,n}, \pi _{j})$, where
$\otimes = \otimes _{L _{n}}$ and $\lambda _{j-1,n}$ is the unique
$L _{n}$-automorphism of $L _{j-1}L _{n}$ extending $\lambda
_{j-1}$, for each index $j$. Therefore, using (3.1) and Lemma
\ref{lemm4.2} (a), one obtains inductively that it suffices to prove
that $T _{n} \in d(K)$, provided $C _{n} \in d(L _{n})$.
\par
Denote by $w _{n}$ the valuation of $C _{n}$ extending $v _{L
_{n}}$, and by $\widehat C _{n}$ its residue division ring. It
follows from the Ostrowski-Draxl theorem that $w _{n}(C _{n})$
equals the sum of $v(M _{n})$ and the group generated by $q
^{-1}v(\pi _{i'})$, $i ^{\prime } = 2, \dots , n - 1$. Similarly, it
is proved that $\widehat C _{n}$ is a field and $\widehat C _{n}
^{q} \subseteq \widehat K$. One also sees that $\widehat C _{n} \neq
\widehat K$ if and only if $\mu < n - 1$, and in this case,
$[\widehat C _{n}\colon \widehat K] = q ^{n-1-\mu }$ and $\hat \pi
_{u} \in \widehat C _{n} ^{q}$, $u = \mu + 1, \dots , n - 1$. These
results show that $v(\pi _{n}) \notin qw _{n}(C _{n})$, if $\mu =
n$, and $\hat \pi _{n} \notin \widehat C _{n} ^{q}$ when $\mu < n$.
Let now $\bar \lambda _{n}$ be the $K$-automorphism of $C _{n}$
extending both $\lambda _{n}$ and the identity of the natural
$K$-isomorphic copy of $T _{n-1}$ in $C _{n}$, and let $t _{n}
^{\prime } = \prod _{\kappa =0} ^{q-1} \bar \lambda _{n} ^{\kappa
}(t _{n})$, for each $t _{n} \in C _{n}$. Then, by Skolem-Noether's
theorem (cf. \cite{P}, Sect. 12.6), $\bar \lambda _{n}$ is induced
by an inner $K$-automorphism of $T _{n}$. This implies $w _{n}(t
_{n}) = w _{n}(\bar \lambda _{n}(t _{n}))$ and $w _{n}(t _{n}
^{\prime }) \in qw _{n}(C _{n})$, for all $t _{n} \in C _{n}$, and
yields $\hat t _{n} ^{\prime } \in \widehat C _{n} ^{q}$ when $w
_{n}(t _{n}) = 0$. Therefore, $t _{n} ^{\prime } \neq \pi _{n}$, $t
_{n} \in C _{n}$, so it follows from \cite{A1}, Ch. XI, Theorems~11
and 12, that $T _{n} \in d(K)$. Lemma \ref{lemm4.2} is proved.
\end{proof}

\medskip
{\it Proof of the latter assertion of Theorem \ref{theo2.1} (c).}
Assume that $F/E$ is an FG-extension, such that char$(E) = p$,
$[E\colon E ^{p}] = p ^{\nu } < \infty $ and trd$(F/E) = t \ge 1$.
This implies $[F\colon F ^{p}] = p ^{\nu +t}$, so it follows from
Lemma \ref{lemm4.1} and \cite{A1}, Ch. VII, Theorem~28, that
Brd$_{p}(F) \le {\rm abrd}_{p}(F) \le \nu + t$. At the same time, it
is clear from (4.1) and Lemma \ref{lemm4.2} that there exists
$\Delta \in d(F)$ with exp$(\Delta ) = p$ and ind$(\Delta ) = p
^{\nu +t-1}$, which yields Brd$_{p}(F) \ge \nu + t - 1$ and so
completes our proof.

\medskip
Our next lemma is implied by (3.5), Lemma \ref{lemm3.1} and the
immediacy of Henselizations of valued fields (cf. \cite{Efr2},
Theorems~15.2.2 and 15.3.5).

\medskip
\begin{lemm}
\label{lemm4.3}
Let $E$ be a field, $F = E(X)$ a rational extension of $E$ with {\rm
trd}$(F/E) = 1$, $f(X) \in E[X]$ an irreducible polynomial over $E$,
$M$ an extension of $E$ generated by a root of $f$ in $E _{\rm sep}$,
$v$ a discrete $E$-valuation of $F$ with a uniform element $f$, and
$(F _{v}, \bar v)$ a Henselization of $(F, v)$. Also, let
$\widetilde D \in d(M)$ be an algebra of exponent $p \in \mathbb P$.
Then $M$ is $E$-isomorphic to the residue field of $(F, v)$ and $(F
_{v}, \bar v)$, and there exists $D \in d(F)$ with {\rm exp}$(D) =
p$ and $[D \otimes _{F} F _{v}] = [D ^{\prime }]$, where $D ^{\prime
} \in d(F _{v})$ is an inertial lift of $\widetilde D$ over $F
_{v}$.
\end{lemm}

\medskip
{\it Proof of Theorem \ref{theo2.1} (a).} Let abrd$_{p}(E) = \lambda
\in \mathbb N$ and $F = E(X _{1}, \dots , X _{\kappa })$. Then, by
Lemma \ref{lemm4.1}, there exists $M \in {\rm Fe}(E)$, such that
$d(M)$ contains an algebra $\widetilde \Delta $ with exp$(\widetilde
\Delta ) = p$ and ind$(\widetilde \Delta ) = p ^{\lambda }$. We show
that there is $\Delta \in d(F)$ with exp$(\Delta ) = p$ and
ind$(\Delta ) \ge p ^{\lambda + \kappa -1}$. Suppose first that
$\kappa = 1$, take a primitive element $\alpha $ of $M/E$, and
denote by $f(X _{1})$ its minimal monic polynomial over $E$. Attach
to $f$ a discrete valuation $v$ of $F$ and fix $(F _{v}, \bar v)$ as
in Lemma \ref{lemm4.3}. Then, by Lemma \ref{lemm3.1}, there is
$\Delta _{1} \in d(F)$ with $[\Delta _{1} \otimes _{F} F _{v}] =
[\overline {\Delta }]$ and exp$(\Delta _{1}) = p$, where $\overline 
{\Delta }$ is an inertial lift of $\widetilde {\Delta }$ over $F
_{v}$. Since $\overline \Delta \in d(F _{v})$ and ind$(\overline
\Delta ) = p ^{\lambda }$, this indicates that $p ^{\lambda } \mid
{\rm ind}(\Delta _{1})$, which proves Theorem \ref{theo2.1} (a) when 
$\kappa = 1$. In addition, Lemma \ref{lemm3.2} implies that there
exist infinitely many degree $p$ cyclic extensions of $F$ in $F
_{v}$. Hence, $F _{v}$ contains as a subfield a Galois extension $R 
_{\kappa }$ of $F$ with $\mathcal{G}(R _{\kappa }/F)$ of order $p 
^{\kappa-1}$ and period $p$. When ind$(\Delta _{1}) = p ^{\lambda
}$, this makes it easy to deduce the existence of $\Delta $, for an 
arbitrary $\kappa $, from (4.1) (with a ground field $E(X _{1})$
instead of $E$) and \cite{Mo}, Theorem~1, or else, by repeatedly
using the Proposition in \cite{P}, Sect. 19.6. It remains to consider
the case where $\kappa \ge 2$ and there exists $D _{1} \in d(E(X
_{1}))$ with exp$(D _{1}) = p$ and ind$(D _{1}) = p ^{\lambda '} > p 
^{\lambda }$. It is easily verified that $D _{1} \otimes _{E(X _{1})}
E(X _{1}) ((X _{2})) \in d(E(X _{1}) ((X _{2})))$, and it follows
from Lemma \ref{lemm3.2} that there are infinitely many degree $p$
cyclic extensions of $E(X _{1}, X _{2})$ in $E(X _{1}) ((X _{2}))$.
As in the case of $\kappa = 1$, this enables one to prove the
existence of $\Delta ^{\prime } \in d(F)$ with exp$(\Delta ^{\prime
}) = p$ and ind$(\Delta ^{\prime }) = p ^{\lambda '+\kappa -2} \ge p 
^{\lambda +\kappa -1}$. Thus Theorem \ref{theo2.1} (a) is proved.

\medskip
\begin{coro}
\label{coro4.4}
Let $E$ be a field and $F/E$ a rational
extension with {\rm trd}$(F/E)$ $= \infty $. Then {\rm Brd}$_{p}(F)
= \infty $, for every $p \in \mathbb P$.
\end{coro}

\medskip
\begin{proof}
This follows from Theorem \ref{theo2.1} (a) and the fact that, for
any rational field extension $F ^{\prime }/F$ with trd$(F ^{\prime
}/F) = 2$, there is an $E$-isomorphism $F \cong F ^{\prime }$,
whence Brd$_{p}(F) = {\rm Brd}_{p}(F ^{\prime })$, for each $p \in
\mathbb P$.
\end{proof}

\medskip
\begin{rema}
\label{rema4.5} Let $E$ be a field with abrd$_{p}(E) = \infty $, $p
\in \mathbb P$, and let $F/E$ be a transcendental FG-extension. Then
it follows from (1.1) (b), (c) and Theorem \ref{theo2.1} (b) that
Brauer pairs $(m, n) \in \mathbb N ^{2}$ are index-exponent pairs
over $F$. Therefore, Corollary \ref{coro4.4} with its proof implies
the latter assertion of (1.2).
\par
Alternatively, it follows from Galois theory, Lemmas \ref{lemm3.2},
\ref{lemm4.3} and basic theory of valuation prolongations that $r
_{p}(\Phi ) = \infty $, $p \in \mathbb P$, for every transcendental
FG-extension $\Phi /E$. Hence, by \cite{Du} and Witt's lemma (cf.
\cite{Dr1}, Sect. 15, Lemma~2), finite abelian groups are realizable
as Galois groups over $\Phi $, so both parts of (1.2) can be proved
by the method used in \cite{P}, Sect. 19.6.
\end{rema}

\medskip
\begin{prop}
\label{prop4.6} Let $F/E$ be an {\rm FG}-extension with {\rm
trd}$(F/E) = t \ge 1$ and {\rm abrd}$_{p}(E) < \infty $, $p \in P$,
for some subset $P \subseteq \mathbb P$. Then $P$ possesses a finite
subset $P(F/E)$, such that {\rm Brd}$_{p}(F) \ge {\rm abrd}_{p}(E) +
t - 1$, $p \in P \setminus P(F/E)$.
\end{prop}

\medskip
{\it Proof.} It follows from (1.1) (c) and Theorem \ref{theo2.1} (a) 
that one may take as $P(F/E)$ the set of divisors of $[F\colon F 
_{0}]$ lying in $P$, for some rational extension $F _{0}$ of $E$ in 
$F$ with trd$(F _{0}/E) = t$.

\medskip
\begin{exam}
\label{exam4.7} There exist field extensions $F/E$ satisfying the
conditions of Proposition \ref{prop4.6}, for $P = \mathbb P$, such
that $P(F/E)$ is nonempty. For instance, let $E$ be a real closed
field, $\Phi $ the function field of the Brauer-Severi variety
attached to the symbol $E$-algebra $A = A _{-1}(-1, -1; E)$, and
$F/\Phi $ a finite field extension with $\sqrt{-1} \notin F$. Then
abrd$(F) = 0 < {\rm abrd}_{2}(E) = 1$ (see the example in
\cite{Ch3}) and abrd$_{p}(E) = 0$, $p > 2$, which implies $P(F/E) =
\{2\}$ and $P = \mathbb P$.
\end{exam}

\medskip
\section{\bf Proof of Theorem \ref{theo2.1} (b)}

\medskip
The former claim of Theorem \ref{theo2.1} (b) is implied by the
following lemma.

\medskip
\begin{lemm}
\label{lemm5.1}
Let $K$ be a field with {\rm abrd}$_{p}(K) = \infty
$, for some $p \in \mathbb P$, and let $F/K$ be an {\rm
FG}-extension with {\rm trd}$(F/K) \ge 1$. Then there exist $D _{\nu
} \in d(F)$, $\nu \in \mathbb N$, such that {\rm exp}$(D _{\nu }) =
p$ and {\rm ind}$(D _{\nu }) \ge p ^{\nu }$.
\end{lemm}

\medskip
\begin{proof}
Statement (1.1) (c) implies the class of fields $\Phi $ with
abrd$_{p}(\Phi ) = \infty $ is closed under the formation of finite
extensions. Since $K$ has a rational extension $F _{0}$ in $F$ with
trd$(F _{0}/K) = {\rm trd}(F/K)$, whence $[F\colon F _{0}] < \infty
$, this shows that it is sufficient to prove Lemma \ref{lemm5.1} in
the case of $F = F _{0}$. Note also that ind$(T _{0} \otimes _{K} F
_{0}) = {\rm ind}(T _{0})$ and exp$(T _{0} \otimes _{K} F _{0}) =
{\rm exp}(T _{0})$, for each $T _{0} \in d(K)$, so one may assume,
for the proof, that $F = F _{0}$ and trd$(F/K) = 1$. It follows from
Lemma \ref{lemm4.1} and the equality abrd$_{p}(K) = \infty $ that
there are $M _{\nu } \in {\rm Fe}(K)$ and $\widetilde D _{\nu } \in
d(M _{\nu })$, $\nu \in \mathbb N$, with exp$(\widetilde D _{\nu })
= p$ and ind$(\widetilde D _{\nu }) \ge p ^{\nu }$, for each index
$\nu $. Hence, by Lemmas \ref{lemm4.3} and \ref{lemm3.1}, there
exist a discrete $K$-valuation $v _{\nu }$ of $F$, and an algebra $D
_{\nu } \in d(F)$, such that the residue field of $(F, v _{\nu })$
is $K$-isomorphic to $M _{\nu }$, exp$(D _{\nu }) = p$, and $[D
_{\nu } \otimes _{F} F _{v _{\nu }}] = [D _{\nu } ^{\prime }]$, where 
$D _{\nu } ^{\prime }$ is an inertial lift of $\widetilde D _{\nu }$
over $F _{v _{\nu }}$. This implies ind$(\widetilde D _{\nu }) \mid 
{\rm ind}(D _{\nu })$, $\nu \in \mathbb N$, proving Lemma 
\ref{lemm5.1}.
\end{proof}

\medskip
To prove the latter part of Theorem \ref{theo2.1} (b) we need the
following lemma.

\medskip
\begin{lemm}
\label{lemm5.2}
Let $A$, $B$ and $C$ be algebras over a field $F$,
such that $A, B, C \in s(F)$, $A = B \otimes _{F} C$, {\rm exp}$(C)
= p \in \mathbb P$, and {\rm exp}$(B) = {\rm ind}(B) = p ^{m}$, for
some $m \in \mathbb N$. Assume that {\rm ind}$(A) = p ^{n} > p ^{m}$
and $k$ is an integer with $m < k \le n$. Then there exists $T _{k}
\in s(F)$ with {\rm exp}$(T _{k}) = p ^{m}$ and {\rm ind}$(T _{k}) =
p ^{k}$.
\end{lemm}

\medskip
\begin{proof}
When $k = n$, there is nothing to prove, so we assume that $k < n$.
By \cite{M1}, Sect. 4, Theorem~2, $[C] = [\Delta _{1} \otimes _{F}
\dots \otimes _{F} \Delta _{\nu }]$, where $\nu \in \mathbb N$ and
for each index $j$, $\Delta _{j} \in d(F)$ and ind$(\Delta _{j}) =
p$. Put $T _{j} = B \otimes _{F} (\Delta _{1} \otimes _{F} \dots
\otimes _{F} \Delta _{j})$ and $t _{j} = {\rm deg}(T _{j})/{\rm
ind}(T _{j})$, $j = 1, \dots , \nu $, and let $S(A)$ be the set of
those $j$, for which ind$(T _{j}) \ge p ^{k}$. Clearly, $S(A) \neq
\phi $ and the set $S _{0}(A) = \{i \in S(A)\colon \ t _{i} \le t
_{j}, j \in S(A)\}$ contains a minimal index $\gamma $. The
conditions of Lemma \ref{lemm5.2} ensure that exp$(T _{j}) = p
^{m}$, so ind$(T _{j}) = p ^{m(j)}$, where $m(j) \in \mathbb N$, for
each $j \in S(A)$. We show that ind$(T _{\gamma }) = p ^{k}$. If
$\gamma = 1$, then (1.1) (c) and the inequality $m < k$ imply $k = m
+ 1$ and ind$(T _{1}) = p ^{k}$, as claimed. Suppose now that
$\gamma \ge 2$. Then it follows from (1.1) (b) that ind$(T _{\gamma
}) = {\rm ind}(T _{\gamma -1}).p ^{\mu }$, for some $\mu \in \{-1,
0, 1\}$. The possibility that $\mu \neq 1$ is ruled out, since it
contradicts the fact that $\gamma \in S _{0}(A)$. This yields ind$(T
_{\gamma }) = {\rm ind}(T _{\gamma -1}).p$ and $t _{\gamma } = t
_{\gamma -1}$. As $\gamma $ is minimal in $S _{0}(A)$, it is now
easy to see that ind$(T _{\gamma -u}) = p ^{k-u}$, $u = 0, 1$, which
proves Lemma \ref{lemm5.2}.
\end{proof}

\medskip
The conditions of Lemma \ref{lemm5.2} are fulfilled, for each $m \in
\mathbb N$ and infinitely many integers $n > m$, if char$(E) = p$,
$E$ is \emph{not} virtually perfect and $F/E$ satisfies the
conditions of Theorem \ref{theo2.1}. Since, by Witt's lemma, cyclic
$p$-extensions of $F$ are realizable as intermediate fields of
$\mathbb Z _{p}$-extensions of $F$, this can be obtained by applying
(1.1) (b), (4.1) and Lemma \ref{lemm4.2} together with general
properties of cyclic $F$-algebras, see \cite{P}, Sect. 15.1,
Corollary~b and Proposition~b. Thus Theorem \ref{theo2.1} is proved
in the case of $p = {\rm char}(E)$. For the proof of the latter
assertion of Theorem \ref{theo2.1} (b), when $p \neq {\rm char}(E)$,
we need the following lemma.

\medskip
\begin{lemm}
\label{lemm5.3}
Let $K$ be a field and $F/K$ an {\rm FG}-extension
with {\rm trd}$(F/K) = 1$. Then, for each $p \in \mathbb P$
different from char$(K)$, there exist non-equivalent discrete
$K$-valuations $v _{m}$ of $F$, $m \in \mathbb N$, satisfying the
following:
\par
{\rm (a)} For any $m \in \mathbb N$, $(F, v _{m})$ possesses a
totally ramified extension $(F _{m}, w _{m})$, such that $F _{m} \in
I(F _{\rm sep}/F)$, $F _{m}/F$ is cyclic and $[F _{m}\colon F] = p
^{m}$;
\par
{\rm (b)} The valued fields $(F _{m}, w _{m})$ can be chosen so that
$F _{m'} \cap F _{\bar m} = F$, $m ^{\prime } \neq \bar m$.
\end{lemm}

\medskip
\begin{proof}
Let $X \in F$ be a transcendental element over $K$. Then $F/K(X)$ is
a finite extension, and the separable closure of $K(X)$ in $F$ is
unramified relative to every discrete $K$-valuation of $K(X)$, with
at most finitely many exceptions (up-to an equivalence, see
\cite{CF}, Ch. I, Sect. 5). This reduces the proof of Lemma
\ref{lemm5.3} to the special case of $F = K(X)$. For each $m \in
\mathbb N$, let $\delta _{m} \in F _{\rm sep}$ be a primitive $p
^{m}$-th root of unity, $K _{m} = K(\delta _{m})$, $f _{m}(X) \in
K[X]$ the minimal polynomial of $\delta _{m}$ over $K$, and $\rho
_{m}$ a discrete $K$-valuation of $F$ with a uniform element $f
_{m}$. Clearly, the valuations $\rho _{m}$, $m \in \mathbb N$, are
pairwise non-equivalent. Also, it is well-known (see \cite{L}, Ch. V,
Theorem 6; Ch. VIII, Sect. 3, and \cite{IR}, Ch. 4, Sect. 1) that if
$m ^{\prime }, \bar m \in \mathbb N$, then the extension $K
_{m'}(\delta _{\bar m})/K _{m'}$ are cyclic except, possibly, in the
case where $m ^{\prime } = 1$, $\bar m > 2$, $p = 2$, char$(K) = 0$
and $\delta _{2} \notin K$. Denote by $v _{m}$ the valuation $\rho
_{m+1}$, for each $m$, if $p = 2$, char$(K) = 0$ and $\delta _{2}
\notin K$, and put $v _{m} = \rho _{m}$, $m \in \mathbb N$,
otherwise. Since $p \neq {\rm char}(K)$, and by Lemma 4.5, $K _{m}$
is $K$-isomorphic to the residue field of $(F, \rho _{m})$, we have
$\delta _{m} \in F _{v _{m}}$, where $F _{v _{m}}$ is a
Henselization of $F$ in $F _{\rm sep}$ relative to $v _{m}$. This
enables one to deduce from Kummer theory that $F _{v _{m}}$
possesses a totally ramified cyclic extension $L _{v _{m}}$ of
degree $p ^{m}$. Furthermore, it follows from the choice of $v _{m}$
and the observation on the extensions $K _{m'}(\delta _{\bar m})/K
_{m'}$ that $F _{v _{m'}}(\delta _{\bar m})/F _{v _{m'}}$ are cyclic, 
for all pairs $m ^{\prime }, \bar m \in \mathbb N$. Hence, by the
generalized Grunwald-Wang theorem (cf. \cite{LR}, Theorems~1 (ii)
and 2) and the note preceding the statement of Lemma \ref{lemm3.2},
there exist totally ramified extensions $(F _{m}, w _{m})/(F, v
_{m})$, $m \in \mathbb N$, such that $F _{m} \in I(F _{\rm sep}/F)$,
$F _{m}/F$ is cyclic with $[F _{m}\colon F] = p ^{m}$, for each $m$,
and in case $m \ge 2$, $F _{m}/F$ is unramified relative to $v _{1},
\dots , v _{m-1}$. This ensures that $F _{m'} \cap F _{\bar m} = F$,
$m ^{\prime } \neq \bar m$, and so completes the proof of Lemma
\ref{lemm5.3}.
\end{proof}

\medskip
{\it Proof of the latter statement of Theorem \ref{theo2.1} (b).}
Let abrd$_{p}(E) = \infty $, for some $p \in \mathbb P$. In view of
(1.1) (b), Lemmas \ref{lemm3.1}, \ref{lemm5.1} and \ref{lemm5.2}, it
is sufficient to show that there exists $A _{m} \in d(F)$ with
exp$(A _{m}) = {\rm ind}(A _{m}) = p ^{m}$, for any fixed $m \in
\mathbb N$. As in the proof of Lemma \ref{lemm5.1}, our
considerations reduce to the special case of trd$(F/K) = 1$.
Analyzing this proof, one obtains that there is $M \in {\rm Fe}(E)$,
such that $d(M)$ contains a cyclic $M$-algebra $\widetilde A _{1}$
of degree $p$, and when $p \neq {\rm char}(E)$, $M$ contains a
primitive $p ^{m}$-th root of unity $\delta _{m}$. Note further that
$M$ can be chosen so as to be $E$-isomorphic to the residue field
$\widehat F$ of $F$ relative to some discrete $E$-valuation $v$. In
view of Kummer theory (see \cite{L}, Ch. VIII, Sect. 6) and Witt's
lemma, the assumptions on $M$ ensure that each degree $p$ cyclic
extension $Y _{1}$ of $M$ lies in $I(Y _{m}/M)$, for some degree $p
^{m}$ cyclic extension $Y _{m}/M$. Suppose now that $Y _{1}$ embeds
in $\widetilde A _{1}$ as an $M$-subalgebra, fix a generator $\tau
_{1}$ of $\mathcal{G}(Y _{1}/M)$ and an automorphism $\tau _{m}$ of
$Y _{m}$ extending $\tau _{1}$. Then $\widetilde A _{1}$ is
isomorphic to the cyclic $M$-algebra $(Y _{1}/M, \tau _{1}, \tilde
\beta )$, for some $\tilde \beta \in M ^{\ast }$, $\tau _{m}$
generates $\mathcal{G}(Y _{m}/M)$, the $M$-algebra $\widetilde A
_{m} = (Y _{m}/M, \tau _{m}, \tilde \beta )$ lies in $s(M)$, and we
have $p ^{m-1}[\widetilde A _{m}] = [\widetilde A _{1}]$ (cf.
\cite{P}, Sect. 15.1, Corollary~b). Therefore, $\widetilde A _{m}
\in d(M)$ and ind$(\widetilde A _{m}) = {\rm exp}(\widetilde A _{m})
= p ^{m}$. Assume now that $(F, v)$ has a valued extension $(L, v
_{L})$, such that $L/F$ is cyclic, $[L\colon F] = p ^{m}$ and the
residue field of $(L, v _{L})$ is $E$-isomorphic to $Y _{m}$. Then
$\mathcal{G}(L/F) \cong \mathcal{G}(Y _{m}/M)$, and for each
generator $\sigma $ of $\mathcal{G}(L/F)$ and pre-image $\beta $ of
$\tilde \beta $ in $O _{v}(F)$, the algebra $A _{m} = (L/F, \sigma ,
\beta )$ lies in $d(F)$ (see \cite{P}, Sect. 15.1, Proposition~b,
and \cite{JW}, Theorem~5.6). Note also that ind$(A _{m}) = {\rm
exp}(A _{m}) = p ^{m}$ and $\sigma $ can be chosen so that $A _{m}
\otimes _{F} F _{v}$ be an inertial lift of $\widetilde A _{m}$ over
$F _{v}$. When $p > 2$, this completes the proof of Theorem
\ref{theo2.1} (b), since Lemma \ref{lemm3.2} guarantees in this case
the existence of a valued extension $(L, v _{L})$ of $(F, v)$ with
the above-noted properties.
\par
Similarly, one concludes that if $p = 2$, then it suffices to prove
Theorem \ref{theo2.1} (b), provided char$(E) = 0$ and 
$\mathcal{G}(E(\delta _{m})/E)$ is noncyclic, where $\delta _{m}$ is 
a primitive $2 ^{m}$-th root of unity in $E _{\rm sep}$. This 
implies the group $E _{1} ^{\ast }/E _{1} ^{\ast 2 ^{\nu }}$ has
period $2 ^{\nu }$, for each $\nu \in \mathbb N$, $E _{1} \in {\rm
Fe}(E)$ (cf. \cite{L}, Ch. VIII, Sects. 3 and 9). Take a valued
extension $(F _{m}, w _{m})/(F, v _{m})$ as required by Lemma
\ref{lemm5.3}, and denote by $\widehat F _{m}$ the residue field of
$(F, v _{m})$. Fix a generator $\psi _{m}$ of $\mathcal{G}(F
_{m}/F)$ and an element $\tilde \beta _{m} \in \widehat F _{m}
^{\ast }$ so that $\tilde \beta _{m} ^{2 ^{m-1}} \notin \widehat F
_{m} ^{\ast 2 ^{m}}$, and put $A _{m} = (F _{m}/F, \psi _{m}, \beta
_{m})$, for some pre-image $\beta _{m}$ of $\tilde \beta _{m}$ in $O
_{v _{m}}(F)$. As $(F _{m}, w _{m})/(F, v _{m})$ is totally
ramified, $w _{m}$ is uniquely determined by $v _{m}$, up-to an
equivalence. Therefore, $w _{m}(\lambda _{m}) = w _{m}(\psi
_{m}(\lambda _{m}))$, for all $\lambda _{m} \in F _{m}$, and when $w
_{m}(\lambda _{m}) = 0$, $\widehat F _{m} ^{\ast 2 ^{m}}$ contains
the residue class of the norm $N _{F} ^{F _{m}}(\lambda _{m})$. Now
it follows from \cite{P}, Sect. 15.1, Proposition~b, that $A _{m}
\in d(F)$ and ind$(A _{m}) = {\rm exp}(A _{m}) = 2 ^{m}$, so Theorem
\ref{theo2.1} is proved.

\medskip
\begin{coro}
\label{coro5.4}
Let $E$ be a field with {\rm abrd}$(E) = \infty
$. Then {\rm Brd}$(F) = \infty $, for every transcendental {\rm
FG}-extension $F/E$.
\end{coro}

\medskip
\begin{proof}
The equality abrd$(E) = \infty $ means that either abrd$_{p'}(E) =
\infty $, for some $p ^{\prime } \in \mathbb P$, or abrd$_{p}(E)$,
$p \in \mathbb P$, is an unbounded number sequence. In view of
Theorem \ref{theo2.1} (b) and Proposition \ref{prop4.6}, this proves
our assertion.
\end{proof}

\medskip
Corollary \ref{coro5.4} shows that a field $E$ satisfies abrd$(E) <
\infty $, if its FG-extensions are of finite dimensions, in the sense 
of \cite{ABGV}, Sect. 4. In view of (2.7) (a), this proves that
Problem~4.4 of \cite{ABGV} is solved, generally, in the negative,
even when all finite extensions of $E$ have finite Brauer 
dimensions. Statements (2.7) also imply that both cases pointed out 
in the proof of Corollary \ref{coro5.4} can be realized.

\medskip
\begin{rema}
\label{rema5.5} Statement (2.6) indicates that if \cite{ABGV},
Problem~4.5, is solved affirmatively in the class $\mathcal{A}$ of
virtually perfect fields $E$ with abrd$(E) < \infty $, then abrd$(E)
\le {\rm dim}(E)$. We show that such a solvability would imply the
numbers $c(E)$, in (2.6), depend on the choice of $E$ and may be
arbitrarily large. Let $C$ be an algebraically closed field, $\nu $
a positive integer and $C _{\nu } = C((X _{1})) \dots ((X _{\nu }))$
the iterated formal Laurent formal power series field in $\nu $
variables over $C$. We prove that $c(C _{\nu }) \ge [\nu /2]
- 1$. Note first that each FG-extension $F/C _{\nu }$ with trd$(F/C
_{\nu }) = 1$ has a $C$-valuation $f _{\nu }$, such that
trd$(\widehat F/C) = 1$ and $f _{\nu }(F) = \mathbb Z ^{\nu }$.
Indeed, if $T \in F$ is a transcendental element over $C _{\nu }$,
$F _{0} = C _{\nu }(T)$, and $f _{0}$ is the restricted Gauss
valuation of $F _{0}$ extending the natural $\mathbb Z ^{\nu
}$-valued $C$-valuation of $C _{\nu }$ (see \cite{Efr2},
Example~4.3.2), then one may take as $f _{\nu }$ any prolongation of
$f _{0}$ on $F$. The equality trd$(\widehat F/C) = 1$ ensures that
$r _{p}(\widehat F) = \infty $, for all $p \in \mathbb P$, which
enables one to deduce from \cite{Mo}, Theorem~1, and \cite{LiKr},
Corollary~1.4, that Brd$_{p}(F) = {\rm abrd}_{p}(F) = \nu $, $p \in
\mathbb P$ and $p \neq {\rm char}(C)$ (see \cite{LiKr}, page 37, for
more details in case $F/C _{\nu }$ is rational). At the same time,
it follows from \cite{Ch5}, Proposition~7.1, that if char$(C) = 0$,
then Brd$(C _{\nu }) = {\rm abrd}(C _{\nu }) = [\nu /2]$; hence, by
(2.6), $c(C _{\nu }) \ge {\rm abrd}(F) - {\rm abrd}(C _{\nu }) - 1 = 
\nu - [\nu /2] - 1 \ge [\nu /2] - 1$, as claimed.
\end{rema}

\medskip
\begin{coro}
Let $F$ be a rational extension of an algebraically closed field $F
_{0}$. Then {\rm trd}$(F/F _{0}) = \infty $ if and only if each
Brauer pair $(m, n) \in \mathbb N ^{2}$ is realizable as an
index-exponent pair over $F$.
\end{coro}

\medskip
\begin{proof}
If trd$(F/F _{0}) = n < \infty $, then finite extensions of $F$ are
$C _{n}$-fields, by Lang-Tsen's theorem \cite{L1}, so Lemma
\ref{lemm4.1} and \cite{Mat} imply Brd$_{p}(F) < p ^{n-1}$, $p \in
\mathbb P$ (see \cite{MS}, (16.10), for case $p = 2$). In view of
(1.2), this completes our proof.
\end{proof}

\medskip
Theorem \ref{theo2.1} and Example \ref{exam4.7} lead naturally to
the question of whether Brd$_{p}(F) \ge k + {\rm trd}(F/E)$,
provided that $F/E$ is an FG-extension and Brd$_{p}(E ^{\prime }) =
k < \infty $, $E ^{\prime } \in {\rm Fe}(E)$, for a given $p \in
\mathbb P$. Our next result gives an affirmative answer to this
question in several frequently used special cases:

\medskip
\begin{prop}
\label{prop5.6}
Let $E$ be a field and $F$ an {\rm
FG}-extension of $E$ with {\rm trd}$(F/E) = n > 0$. Suppose that
there exists $M \in {\rm Fe}(E)$ satisfying the following
condition, for some $p \in \mathbb P$ and $k \in \mathbb N$:
\par
{\rm (c)} For each $M ^{\prime } \in {\rm Fe}(M)$, there are $D
^{\prime } \in d(M ^{\prime })$ and $L ^{\prime } \in I(M ^{\prime
}(p)/M ^{\prime })$, such that {\rm exp}$(D ^{\prime }) = [L
^{\prime }\colon M ^{\prime }] = p$, {\rm ind}$(D ^{\prime }) = p
^{k}$ and $D ^{\prime } \otimes _{M'} L ^{\prime } \in d(L ^{\prime
})$.
\par
Then there exist $D \in d(F)$, such that {\rm exp}$(D) = p$ and {\rm
ind}$(D) \ge p ^{k+n}$; in particular, {\rm Brd}$_{p}(F) \ge k + n$.
\end{prop}

\medskip
Proposition \ref{prop5.6} is proved along the lines drawn in the
proofs of Theorem \ref{theo2.1} (a) and (b), so we omit the details.
Note only that if $n \ge 2$ or $k = 1$, then $D$ can be chosen so
that $D \otimes _{F} F _{v} \in d(F _{v})$, $[D \otimes _{F} F _{v}]
\in {\rm Br}(F _{v,{\rm un}}/F _{v})$ and $p ^{n-1} \mid e(D \otimes
_{F} F _{v}/F _{v}) \mid p ^{n}$, for some $E$-valuation $v$ of $F$
with $\mathbb Z ^{n-1} \le v(F) \le \mathbb Z ^{n}$.

\medskip
\begin{rema}
\label{rema5.8} Condition (c) of Proposition \ref{prop5.6} is
fulfilled, for $k = 1 = {\rm abrd}(E)$ and any $p \in \mathbb P$, if
$E$ is a global field or an FG-extension of an algebraically closed
field $E _{0} ^{\prime }$ with trd$(E/E _{0} ^{\prime }) = 2$. It
also holds when $k = 1$, $p \in \mathbb P$ and $E$ is an
FG-extension of a perfect PAC-field $E _{0}$ with trd$(E/E _{0}) = 1
= {\rm cd}_{p}(E _{0})$ (see \cite{Efr1}, Sect. 3, and \cite{P}, 
Sect. 19.3). In these cases, it can be deduced from (3.1) and 
\cite{Mo}, Theorem~1, that the power series fields $E _{m} = E((X 
_{1})) \dots ((X _{m}))$, $m \in \mathbb N$, satisfy (c), for $k = 1 
+ m = {\rm abrd}_{p}(E _{m})$ (cf. \cite{LiKr}, Appendix~A, or 
\cite{Ch5}, (5.2) and Proposition~5.1). In addition, the conclusion 
of Proposition \ref{prop5.6} is valid, if $E$ is a local field, $k = 
1$ and $p \in \mathbb P$, although (c) is then violated, for every 
$p$ (see Proposition \ref{prop6.3} with its proof, and appendices to 
\cite{Sal3} and \cite{CF}, Ch. VI, Sect. 1).
\end{rema}

\medskip
For a proof of the concluding result of this section, we refer the
reader to \cite{Ch3}. When $F/E$ is a rational extension and $r
_{p}(E) \ge {\rm trd}(F/E)$, this result is contained in \cite{Nak}.
Combined with Lemma \ref{lemm3.2}, it implies Nakayama's 
inequalities
\par\noindent
Brd$_{p'}(F ^{\prime }) \ge {\rm trd}(F ^{\prime }/E ^{\prime }) -
1$, $p ^{\prime } \in \mathbb P$, for any FG-extension $F ^{\prime
}/E ^{\prime }$.

\medskip
\begin{prop}
\label{prop5.9}
Let $F/E$ be an {\rm FG}-extension with {\rm
trd}$(F/E) = n \ge 1$ and {\rm cd}$_{p}(\mathcal{G}_{E}) \neq 0$,
for some $p \in \mathbb P$. Then {\rm Brd}$_{p}(F) \ge n$ except,
possibly, if $p = 2$, the Sylow pro-$2$-subgroups of
$\mathcal{G}_{E}$ are of order $2$, and $F$ is a nonreal field.
\end{prop}

\medskip
It is not known whether an FG-extension $F/E$ with trd$(F/E) = n \ge
3$ satisfies abrd$_{p}(F) = {\rm Brd}_{p}(F) = n - 1$, provided that  
$p \in \mathbb P$, cd$_{p}(\mathcal{G}_{E}) = 0$, and $E$ is perfect
in the case where $p = {\rm char}(E)$. It follows from (1.1) (c)
that this question is equivalent to the Standard Conjecture on $F/E$
(stated by Colliot-Th\'{e}l\`{e}ne, see \cite{LiKr} and \cite{Lieb},
Sect. 1) when $E$ is algebraically closed. The question is also open
in the case excluded by Proposition \ref{prop5.9}. Results like
\cite{Mat}, Theorem~6.3 and Corollary~7.3, as well as statements
(2.1) and (2.3) attract interest in the problem of finding exact
upper bounds on abrd$_{p}(F)$, $p \in \mathbb P$. Specifically, it
is worth noting that if $E$ is algebraically closed and Brd$_{p}(F)
\ge p ^{n-2}$, for infinitely many $p \in \mathbb P$, then this would
solve negatively \cite{ABGV}, Problem~4.5, by showing that Brd$(F) =
\infty $ whenever $n \ge 3$.

\medskip
\section{\bf Reduction of (2.2) to the case of char$(E) = 0$}

\medskip
In this section we show that if $\mathcal{C}$ is a class of
profinite groups and $n$ is a positive integer, then the answer to
(2.2) would be affirmative, for FG-extensions $F/E$ with
$\mathcal{G}_{E} \in \mathcal{C}$ and trd$(F/E) \le n$, if this
holds when char$(E) = 0$. This result can be viewed as a refinement
of \cite{Efr2}, Corollary~22.2.3, in the spirit of \cite{Lieb},
4.1.2.

\medskip
\begin{prop}
\label{prop6.1} Let $E$ be a field of characteristic $q > 0$ and
$F/E$ an {\rm FG}-extension. Then there exists an {\rm FG}-extension
$L/E ^{\prime }$ satisfying the following:
\par
{\rm (a)} {\rm char}$(E ^{\prime }) = 0$, $\mathcal{G}_{E'} \cong
\mathcal{G}_{E}$ and {\rm trd}$(L/E ^{\prime }) = {\rm trd}(F/E)$;
\par
{\rm (b)} {\rm Brd}$_{p}(L) \ge {\rm Brd}_{p}(F)$, {\rm
abrd}$_{p}(L) \ge {\rm abrd}_{p}(F)$, {\rm Brd}$_{p}(E ^{\prime }) =
{\rm Brd}_{p}(E)$ and {\rm abrd}$_{p}(E ^{\prime }) = {\rm
abrd}_{p}(E)$, for each $p \in \mathbb P$ different from $q$.
\end{prop}

\medskip
\begin{proof}
Fix an algebraic closure $\overline F$ of $F$ and denote by $E _{\rm
ins}$ the perfect closure of $E$ in $\overline F$. The extension $E
_{\rm ins}/E$ is purely inseparable, so it follows from the
Albert-Hochschild theorem (cf. \cite{Se}, Ch. II, 2.2) that the
scalar extension map of Br$(E)$ into Br$(E _{\rm ins})$ is
surjective. Since finite extensions of $E$ in $E _{\rm ins}$ are of
$q$-primary degrees, one obtains from (1.1) (c) that ind$(D \otimes
_{E} E _{\rm ins}) = {\rm ind}(D)$ and exp$(D \otimes _{E} E _{\rm
ins}) = {\rm exp}(D)$, provided $D \in d(E)$ and $q \dagger {\rm
ind}(D)$. Therefore, Brd$_{p}(E) = {\rm Brd}_{p}(E _{\rm ins})$ and
abrd$_{p}(E) = {\rm abrd}_{p}(E _{\rm ins})$, for each $p \in
\mathbb P$, $p \neq q$. As $\mathcal{G}_{E _{\rm ins}} \cong
\mathcal{G}_{E}$ (see \cite{L}, Ch. VII, Proposition~12) and $FE
_{\rm ins}/E _{\rm ins}$ is an FG-extension, this reduces the proof
of Proposition \ref{prop6.1} to the case where $E$ is perfect. It is
known (cf. \cite{Efr2}, Theorems~12.4.1 and 12.4.2) that then there
exists a Henselian field $(K, v)$ with char$(K) = 0$ and $\widehat K
\cong E$, which can be chosen so that $v(K) = \mathbb Z$ and $v(q) =
1$. Moreover, it follows from (3.4), \cite{MT} and Galois theory
(see also the proof of \cite{Efr2}, Corollary~22.2.3) that there is
$E ^{\prime } \in I(K _{\rm sep}/K)$, such that $E ^{\prime } \cap K
_{\rm ur} = K$ and $E ^{\prime }K _{\rm ur} = K _{\rm sep}$. This
ensures that $v(E ^{\prime }) = \mathbb Q$, $\widehat E ^{\prime } =
\widehat K = E$ and $E ^{\prime }_{\rm ur} = E ^{\prime }_{\rm sep}
= K _{\rm sep}$. Hence, by (3.3) and (3.5), $\mathcal{G}_{E'} \cong
\mathcal{G}_{E}$, Brd$_{p}(E ^{\prime }) = {\rm Brd}_{p}(E)$ and
abrd$_{p}(E ^{\prime }) = {\rm abrd}_{p}(E)$, $p \in \mathbb P
\setminus \{q\}$. Observe that, since $E$ is perfect, $F/E$ is
separably generated, i.e. there is $F _{0} \in I(F/E)$, such that $F
_{0}/E$ is rational and $F \in {\rm Fe}(F _{0})$ (cf. \cite{L}, Ch.
X). Note further that each rational extension $L _{0}$ of $E
^{\prime }$ with trd$(L _{0}/E ^{\prime }) = {\rm trd}(F _{0}/E)$
has a restricted Gauss valuation $\omega _{0}$ extending $v _{E'}$
with $\widehat L _{0} = F _{0}$ (cf. \cite{Efr2}, Example~4.3.2).
Fixing $(L _{0}, \omega _{0})$, one can take its valued extension
$(L, \omega )$ so that $L _{\omega } \cong L \otimes _{L _{0}} L
_{0,\omega _{0}}$ is an inertial lift of $F$ over $L _{0,\omega
_{0}}$. This yields $\omega (L) = \omega _{0}(L _{0}) = \mathbb Q$,
$\widehat L \cong F$ over $F _{0}$, $[L\colon L _{0}] = [F\colon F
_{0}]$ and trd$(L/K) = {\rm trd}(F/E)$. It also becomes clear that,
for each $F ^{\prime } \in {\rm Fe}(F)$, there exists a valued
extension $(L ^{\prime }, \omega ^{\prime })$ of $(L, \omega )$ with
$[L ^{\prime }\colon L] = [F ^{\prime }\colon F]$ and $\widehat L
^{\prime } \cong F ^{\prime }$. Observing now that $L ^{\prime }/E
^{\prime }$, $F ^{\prime } \in {\rm Fe}(F)$, are FG-extensions,
applying (3.3) and (3.5) to a Henselization $L ^{\prime } _{\omega
'}$, for any admissible $F ^{\prime }$, and using Lemmas
\ref{lemm3.1} and \ref{lemm4.1}, one concludes that Brd$_{p}(L
^{\prime }) \ge {\rm Brd}_{p}(F ^{\prime })$ and abrd$_{p}(L) \ge
{\rm abrd}_{p}(F)$, for all $p \in \mathbb P \setminus \{q\}$.
Proposition \ref{prop6.1} is proved.
\end{proof}

\medskip
We show that in zero characteristic Proposition \ref{prop2.2} can be
deduced from Proposition \ref{prop6.1}.

\medskip
\begin{exam}
\label{exam6.2}
Let $K _{0}$ be a field with $2$ elements, $K _{n} = K _{0}((X _{1}))
\dots ((X _{n}))$, $n \in \mathbb N$, a sequence of iterated formal 
power series fields in $n$ variables over $K _{0}$, inductively 
defined by the rule $K _{n} = K _{n-1}((X _{n}))$, for each $n \in
\mathbb N$, and let $\Theta $ be a perfect closure of the union $K
_{\infty } = \cup _{n=1} ^{\infty } K _{n}$. It is known that the
natural $\mathbb Z ^{n}$-valued valuations, say $v _{n}$, of the
fields $K _{n}$, $n \in \mathbb N$, extend uniquely to a Henselian
$K _{0}$-valuation $v$ of $K _{\infty }$ with $\widehat K _{\infty }
= K _{0}$ and $v(K _{\infty }) = \cup _{n=1} ^{\infty } v _{n}(K
_{n})$. Since $r _{p}(K _{0}) = 1$, $p \in \mathbb P$, and finite
extensions of $K _{\infty }$ in $\Theta $ are totally ramified and
of $2$-primary degrees over $K _{\infty }$, one deduces from
\cite{Ch4}, Lemma~4.4, that Brd$_{p}(K _{\infty }) = {\rm
Brd}_{p}(\Theta ) = 1$ and abrd$_{p}(K _{\infty }) = {\rm
abrd}_{p}(\Theta ) = \infty $, for every $p > 2$. At the same time,
it follows from Lemma \ref{lemm4.2} that $r _{2}(\Theta ) = \infty 
$. Hence, by Proposition \ref{prop6.1}, there is a field $\Theta 
^{\prime }$ with char$(\Theta ) = 0$, abrd$_{2}(\Theta ^{\prime }) = 
0$, $r _{2}(\Theta ^{\prime }) = \infty $, and Brd$_{p}(\Theta 
^{\prime }) = 1$, abrd$_{p}(\Theta ^{\prime }) = \infty $, $p > 2$. 
Moreover, by the proof of Proposition \ref{prop6.1}, $\Theta 
^{\prime }$ can be chosen so that its roots of unity form a 
multiplicative $2$-group. Put $\Theta _{0} = \Theta ^{\prime }$, 
$\Theta _{k} = \Theta _{k-1}((T _{k}))$, $k \in \mathbb N$, and let 
$\theta _{k}$ be the natural (Henselian) $\mathbb Z ^{k}$-valued 
$\Theta _{0}$-valuation of $\Theta _{k}$, for each index $k$. Fix a 
maximal extension $E _{k}$ of $\Theta _{k}$ in $\Theta _{k,{\rm 
sep}}$ with respect to the property that finite extensions of 
$\Theta _{k}$ in $E _{k}$ have odd degrees and are totally ramified 
over $\Theta _{k}$ relative to $\theta _{k}$. This ensures that 
$\widehat E _{k} = \Theta _{0}$, $E _{k}$ does not contain a 
primitive $\mu $-th root of unity, for any odd $\mu > 1$, the group 
$\theta _{k}(E _{k})/2\theta _{k}(E _{k})$ has order $2 ^{k}$, and 
$\theta _{k}(E _{k}) = p\theta _{k}(E _{k})$, for every $p > 2$. 
Therefore, by \cite{Ch4}, Lemma~4.4, Brd$_{2}(E _{k}) = {\rm 
abrd}_{2}(K) = k$, and by (3.5), Brd$_{p}(E _{k}) = 1$ and 
abrd$_{p}(E _{k}) = \infty $, $p > 2$, whence Brd$(E _{k}) = k$.
\end{exam}

\medskip
Similarly to Remark \ref{rema5.5}, the proofs of Proposition
\ref{prop6.1} and our concluding result demonstrate the
applicability of restricted Gauss valuations in finding lower bounds
on Brd$_{p}(F)$, for FG-extensions $F$ of valued fields $E$ with 
abrd$_{p}(E) < \infty $:

\medskip
\begin{prop}
\label{prop6.3} Let $E$ be a local field and $F/E$ an {\rm
FG}-extension. Then {\rm Brd}$_{p}(F)$ $\ge 1 + {\rm trd}(F/E)$, for
every $p \in \mathbb P$.
\end{prop}

\medskip
\begin{proof}
As Brd$_{p}(F) = 1$ when trd$(F/E) = 0$, we assume that trd$(F/E) =
n \ge 1$. We show that, for each $p \in \mathbb P$, there exists $D
_{p} \in d(F)$, such that exp$(D _{p}) = p$, ind$(D _{p}) = p
^{n+1}$ and $D _{p}$ decomposes into a tensor product of cyclic
division $F$-algebras of degree $p$. Let $\omega $ be the standard
discrete valuation of $E$, $\widehat E$ its residue field, and $F
_{0}$ a rational extension of $E$ in $F$ with trd$(F _{0}/E) = n$.
Considering a discrete restricted Gauss valuation of $F _{0}$
extending $\omega $, and its prolongations on $F$, one obtains that
$F$ has a discrete valuation $v$ extending $\omega $, such that
$\widehat F$ is an FG-extension of $\widehat E$ with trd$(\widehat
F/\widehat E) = n$. Hence, by the proof of Proposition
\ref{prop5.9}, given in \cite{Ch3}, there exist $\Delta _{p}
^{\prime } \in d(\widehat F)$ and a degree $p$ cyclic extension $L
_{p} ^{\prime }/\widehat F$, such that $\Delta _{p} ^{\prime }
\otimes _{\widehat F} L _{p} ^{\prime } \in d(L _{p} ^{\prime })$,
exp$(\Delta _{p} ^{\prime }) = p$, ind$(\Delta _{p} ^{\prime }) = p
^{n}$ and $\Delta _{p} ^{\prime }$ is a tensor product of cyclic
division $\widehat F$-algebras of degree $p$. Given a Henselization
$(F _{v}, \bar v)$ of $(F, v)$, Lemma \ref{lemm3.1} implies the
existence of $\Delta _{p} \in d(F)$, such that $\Delta _{p} \otimes
_{F} F _{v} \in d(F _{v})$ is an inertial lift of $\Delta _{p}
^{\prime }$ over $F _{v}$. Also, by Lemma \ref{lemm3.2}, there is a
degree $p$ cyclic extension $L _{p}/F$ with $L _{p} \otimes _{F} F
_{v}$ an inertial lift of $L _{p} ^{\prime }$ over $F _{v}$. Fix a
generator $\sigma $ of $\mathcal{G}(L _{p}/F)$, take a uniform
element $\beta $ of $(F, v)$, and put $D _{p} = \Delta _{p} \otimes
_{F} (L _{p}/F, \sigma , \beta )$. Then it follows from (3.1) and
\cite{Mo}, Theorem~1, that $D _{p} \in d(F)$, exp$(D _{p}) = p$,
ind$(D _{p}) = p ^{n+1}$ and $D _{p} \otimes _{F} F _{v} \in d(F
_{v})$, so Proposition \ref{prop6.3} is proved.
\end{proof}

\medskip
Note finally that if $E$ is a local field, $F/E$ is an FG-extension
and trd$(F/E) = 1$, then Brd$_{p}(F) = 2$, for every $p \in \mathbb
P$. When $p = {\rm char}(E)$, this is implied by Proposition
\ref{prop6.3} and Theorem \ref{theo2.1} (c), and for a proof in the
case of $p \neq {\rm char}(E)$, we refer the reader to \cite{PaSur},
Theorems~1 and 3, \cite{Sal3} and \cite{LiKr}, Corollary~1.4.

\vskip0.4truecm\noindent
\emph{Acknowledgements.}
The concluding part of this research was done during my visit to
Tokai University, Hiratsuka, Japan, in 2012. I would like to thank
my host-professor Junzo Watanabe, the colleagues at the Department
of Mathematics, and Mrs. Yoko Kinoshita and her team for their
genuine hospitality.

\par


\medskip
\begin{thebibliography}{aa}

\bibitem{A1} Albert, A.A.: \emph{Structure of Algebras}. Am. Math.
Soc. Colloq. Publ., vol. XXIV (1939)

\bibitem{ABGV} Auel, A., Brussel, E., Garibaldi, S., Vishne, U.:
\emph{Open problems on central simple algebras}. Transf. Groups 
{\bf 16} (2011), 219-264

\bibitem{CF} Cassels, J.W.S., Fr\" ohlich A. (Eds.): \emph{Algebraic
Number Theory}. Proc. Instruct. Conf., organized by the London Math.
Soc. (a NATO Adv. Study Inst.) with the support of IMU, Univ. of
Sussex, Brighton, 01.9-17.9, 1965, Academic Press, London-New York 
(1967)

\bibitem{Ch1} Chipchakov, I.D.: \emph{The normality of locally finite
associative division algebras over classical fields}. Vestn. Mosk.
Univ., Ser. I (1988), No. 2, 15-17 (Russian: English transl. in:
Mosc. Univ. Math. Bull. {\bf 43} (1988), 2, 18-21)

\bibitem{Ch2} Chipchakov, I.D.: \emph{On the classification of
central division algebras of linearly bounded degree over global
fields and local fields}. J. Algebra {\bf 160} (1993), 342-379

\bibitem{Ch3} Chipchakov, I.D.: \emph{Lower bounds and infinity
criterion for Brauer $p$-dimensions of finitely-generated field
extensions}. C.R. Acad. Buld. Sci. {\bf 66} (2013), 923-932

\bibitem{Ch4} Chipchakov, I.D.: \emph{On the behaviour of Brauer
$p$-dimensions under finitely-generated field extensions}. J. 
Algebra {\bf 428} (2015), 190-204

\bibitem{Ch5} Chipchakov, I.D.: \emph{On Brauer $p$-dimensions and
absolute Brauer $p$-dimensions of Henselian fields}. Preprint,
arXiv:1207.7120v6 [math.RA] (Sept. 02, 2014)

\bibitem{Dr1} Draxl, P.K.: \emph{Skew Fields}. London Math. Soc.
Lecture Notes, vol. 81, Cambridge University Press IX, Cambridge
etc. (1983)

\bibitem{Dr2} Draxl, P.K.: \emph{Ostrowski's theorem for Henselian
valued skew fields}. J. Reine Angew. Math. {\bf 354} (1984),
213-218

\bibitem{Du} Ducos, L.: \emph{R\'{e}alisation r\'{e}guli\`{e}re
explicite des groupes ab\'{e}liens finis comme groupes de Galois}.
J. Number Theory {\bf 74} (1999), 44-55

\bibitem{Efr1} Efrat, I.: \emph{A Hasse principle for function
fields over PAC fields}. Isr. J. Math. {\bf 122} (2001), 43-60

\bibitem{Efr2} Efrat, I.: \emph{Valuations, Orderings, and Milnor
$K$-Theory}. Math. Surveys and Monographs, 124, Providence, RI:
Am. Math. Soc., XIII (2006)

\bibitem{FV} Fesenko, I.B., Vostokov, S.V.: \emph{Local Fields and
Their Extensions}. 2nd ed., Transl. Math. Monographs, 121, Am.
Math. Soc., Providence, RI (2002)

\bibitem{FJ} Fried, M.J., Jarden, M.: \emph{Field Arithmetic}. 2nd
revised and enlarged ed., Ergebnisse der Math. Und ihrer
Grenzgebiete, 3. Folge, Bd. 11, Springer, Berlin (2005)

\bibitem{HHKr} Harbater, D., Hartmann, J., Krashen, D.:
\emph{Applications of patching to quadratic forms and central simple
algebras}. Invent. Math. {\bf 178} (2009), 231-263

\bibitem{IR} Ireland, K., Rosen, M.: \emph{A Classical Introduction 
to Modern Number Theory}. Graduate Texts in Math., vol. 84,
Springer-Verlag, XIII, New York-Heidelberg-Berlin (1982)

\bibitem{JW} Jacob, B., Wadsworth, A.: \emph{Division algebras
over Henselian fields}. J. Algebra {\bf 128} (1990), 126-179

\bibitem{Jong} De Jong, A.J.: \emph{The period-index problem for the
Brauer group of an algebraic surface}. Duke Math. J. {\bf 123}
(2004), 71-94

\bibitem{Ka} Kahn, B.: \emph{Comparison of some field invariants}.
J. Algebra {\bf 232} (2000), 485-492

\bibitem{Kol} Koll\'{a}r, J.: \emph{A conjecture of Ax and
degenerations of Fano varieties}. Isr. J. Math. {\bf 162} (2007),
235-251

\bibitem{L1} Lang, S.: \emph{On quasi algebraic closure}. Ann. Math.
(2) {\bf 55} (1952), 373-390

\bibitem{L} Lang, S.: \emph{Algebra}. Addison-Wesley Publ. Comp.,
Mass. (1965)

\bibitem{Lieb} Lieblich, M.: \emph{Twisted sheaves and the
period-index problem}. Compos. Math. {\bf 144} (2008), 1-31

\bibitem{LiKr} Lieblich, M.: \emph{Period and index in the Brauer
group of an arithmetic surface. With an appendix by D. Krashen}. J.
Reine Angew. Math. {\bf 659} (2011), 1-41

\bibitem{LR} Lorenz, F., Roquette, P.: \emph{The theorem of
Grunwald-Wang in the setting of valuation theory}. F.-V. Kuhlmann
(ed.) et. al., Valuation theory and its applications, vol. II
(Saskatoon, SK, 1999), 175-212, Fields Inst. Commun., {\bf 33},
Am. Math. Soc., Providence, RI, 2003

\bibitem{Mat} Matzri, E.: \emph{Symbol length in the Brauer group of 
a field}. Trans. Am. Math. Soc. (to appear)

\bibitem{MT} Mel'nikov, O.V., Tavgen', O.I.: \emph{The absolute
Galois group of a Henselian field}. Dokl. Akad. Nauk BSSR {\bf 29}
(1985), 581-583

\bibitem{M1} Merkur'ev, A.S.: \emph{Brauer groups of fields}. 
Commun. Algebra {\bf 11} (1983), 2611-2624

\bibitem{MS} Merkur'ev, A.S., Suslin, A.A.: \emph{$K$-cohomology of
Severi-Brauer varieties and norm residue homomomorphisms}. Izv.
Akad. Nauk SSSR {\bf 46} (1982), 1011-1046 (Russian: English transl.
in: Math. USSR Izv. {\bf 21} (1983), 307-340)

\bibitem{Mo} Morandi, P.: \emph{The Henselization of a valued
division algebra}. J. Algebra {\bf 122} (1989), 232-243

\bibitem{Nak} Nakayama, T.: \emph{\"{U}ber die direkte Zerlegung
eines Divisionsalgebra}. Jap. J. Math. {\bf 12} (1935), 65-70

\bibitem{PaSur} Parimala, R., Suresh, V.: \emph{Period-index and
$u$-invariant questions for function fields over complete discretely
valued fields}. Invent. Math. {\bf 197} (2014), No. 1, 215-235 
(Preprint, arXiv:1304.2214v1 [math.RA] (Apr. 08, 2013)

\bibitem{P} Pierce, R.: \emph{Associative Algebras}. Graduate Texts
in Math., vol. 88, Springer-Verlag, XII, New York-Heidelberg-Berlin
(1982)

\bibitem{Re} Reiner, M.: \emph{Maximal Orders}. London Math. Soc.
Monographs, vol. 5, London-New York-San Francisco: Academic Press, a
subsidiary of Harcourt Brace Jovanovich, Publishers (1975)

\bibitem{Sal2} Saltman, D.J.: \emph{Generic algebras}. Brauer groups
in ring theory and algebraic geometry, Proc., Antwerp, 1981, Lect.
Notes in Math. {\bf 917} (1982), 96-117

\bibitem{Sal3} Saltman, D.J.: \emph{Division algebras over $p$-adic
curves}. J. Ramanujan Math. Soc. {\bf 12} (1997), 25-47 (correction
in: ibid. {\bf 13} (1998), 125-129)

\bibitem{Sch} Schilling, O.F.G.: \emph{The Theory of Valuations}.
Mathematical Surveys, No. 4, Am. Math. Soc., New York, N.Y. 
(1950)

\bibitem{Se} Serre, J.-P.: \emph{Galois Cohomology}. Transl. from the
French original by Patrick Ion, Springer, Berlin (1997)

\bibitem{Vo} Voevodsky, V.: \emph{On motivic cohomology with $\mathbb 
Z/l$-coefficients}. Ann. Math. {\bf 174} (2011), 401-438

\end{thebibliography}
\end{document}